\documentclass[a4paper,fleqn]{cas-sc}

\usepackage[numbers]{natbib}

\usepackage{graphicx}
\usepackage{color,soul}
\usepackage{array}
\usepackage{amsmath}
\usepackage{amsfonts}
\usepackage{amssymb}
\usepackage{psfrag}
\usepackage{epstopdf}
\usepackage{subcaption}
\usepackage{booktabs}
\usepackage{latexsym}
\usepackage{kantlipsum}
\usepackage{balance}
\usepackage{lscape}
\usepackage{rotating}
\usepackage{multicol}
\usepackage{multirow,bigdelim}
\usepackage{color, colortbl}
\usepackage{threeparttable}
\usepackage{float}

\usepackage{tikz}
\usetikzlibrary{calc,fit,shapes,arrows,positioning,chains,arrows.meta}
\usepackage{makeshape}
\usepackage{caption}
\usepackage{floatflt}

\newcommand{\parder}[2]{\frac{\partial #1}{\partial #2}}

\newcommand{\eval}[2][\right]{\relax\ifx#1\right\relax \left.\fi#2#1\rvert}

\newcommand{\bv}{\boldsymbol{v}}

\newcommand{\bB}{\boldsymbol{B}}

\newcommand{\bE}{\boldsymbol{E}}

\newcommand{\del}{\boldsymbol{\nabla}}
\newcommand{\bcross}{\boldsymbol{\times}}

\makeatletter
\newcommand{\Rmnum}[1]{\expandafter\@slowromancap\romannumeral #1@}
\makeatother

\newcommand{\bJ}{\boldsymbol{J}}

\newcommand{\bj}{\boldsymbol{j}}

\def\tsc#1{\csdef{#1}{\textsc{\lowercase{#1}}\xspace}}
\tsc{WGM}
\tsc{QE}
\tsc{EP}
\tsc{PMS}
\tsc{BEC}
\tsc{DE}

\begin{document}
\let\WriteBookmarks\relax
\def\floatpagepagefraction{1}
\def\textpagefraction{.001}
\shorttitle{}
\shortauthors{Durgarao Kamireddy, Arup Nandy}

\title [mode = title]{A Novel Conversion Technique from Nodal  to Edge Finite Element Data Structure for Electromagnetic Analysis}
\tnotemark[1]

\tnotetext[1]{Supported by Science \& Engineering Research Board (SERB), and Department of Science \& Technology (DST), Government of India, under the project IMP/2019/000276.}




\author[]{Durgarao Kamireddy}[type=editor,
                        orcid=0000-0002-8546-6306]
\cormark[1]
\ead{durga176103010@iitg.ac.in}


\address[]{Department of Mechanical Engineering, Indian Institute of Technology Guwahati, Guwahati 781039, India}

\author[]{Arup Nandy}
\ead{arupn@iitg.ac.in}

\cortext[cor1]{Corresponding author}


\begin{abstract}
Standard nodal finite elements in electromagnetic analysis have well-known limitation of occurrence of spurious solution. In order to circumvent the problem, a penalty function method or a regularization method is used with potential formulation. These methods solve the problem partially by pushing the spurious mode to the higher end of the spectrum. But it fails to capture singular eigen values in case of the problem domains with sharp edges and corners. To circumvent this limitation, edge elements have been developed for electromagnetic analysis where degree of freedoms are along the edges. But most of the preprocessors develop complex meshes in nodal framework. In this work, we have developed a novel technique to convert nodal data structure to edge data structure for electromagnetic analysis. We have explained the conversion algorithm in details, mentioning associated complexities with relevant examples. The performance of the developed algorithm has been demonstrated extensively with several examples.
\end{abstract}




\maketitle
\section{Introduction}

	The finite element method (FEM) has been widely used for radiation and scattering problems in interior and exterior domains which has large applications in antenna radiations, waveguide transmissions etc. In order to apply the FEM technique the domain can be discretized with either edge element or nodal element. Problem of occurrence of spurious solution is well-known limitation of standard nodal finite elements in electromagnetic analysis. In order to circumvent the problem, the penalty function method and regularization method ~\cite{Paulsen1991},~\cite{Member1985},~\cite{Otin2010} have been used extensively in nodal FEM framework. These methods solve the problem partially by pushing the spurious mode to the higher end of the spectrum. In order to capture singular eigen value in case of sharp edges and corners, in the regularization method we have to use a penalty parameter varying from 0 at the sharp edge to 1 at a large distance from the sharp edge.Also, in order to take care of inherent tangential continuity and normal discontinuity across material interface, potential formulation is used in nodal finite element framework ~\cite{Nandy2016},~\cite{arupthesis},~\cite{Nandy2018},~\cite{Agrawal2017},~\cite{Nandy2018a}. In~\cite{Jog2014}, a two field variation formulation in electromagnetics which can predict the eigen frequencies very accurately with correct multiplicities. There is no ad-hoc term in the mixed formulation as in the penalty function or the regularization method. This method worked very well for all two dimensional geometries like non-convex domains with sharp corners, in-homogeneous domains, curved domains etc. In three dimension, mixed FEM worked for plane structures (structures without any curvature) quite well; there, it worked flawlessly with sharp edges and in-homogeneous domains. But in the case of curved three dimensional geometries, this mixed formulation failed. 
\par      Edge elements were introduced by Whitney, which are also called curl-conforming elements. J C Nedelec presented the conceptual theory of edge elements~\cite{Jinming-finite}. He presented a non conforming tetrahedron and cube finite elements construction~\cite{Nedelec1980} conforming the H curl and H div spaces. Whitney spaces can act as bases for edge elements in FEM for field type of problems~\cite{Bossavit1988} and eddy current problems. In~\cite{Cendes1991},~\cite{Boffi2001},~\cite{Boffi2010},~\cite{Ainsworth2003},~\cite{Garcia-Castillo2000},~\cite{Bramble2005},~\cite{Kamireddy_2021},~\cite{Kamireddy2020} edge elements are used to solve the eigenvalue problems with different shapes. In~\cite{Bossavit1982}, Alain Bossavit et al. solved the 3-D eddy current problems using the combination of FEM and Boundary Integral Element Method (BIEM) methods. Edge elements can be easily applied to an exterior domain problems where the coupling of other methods like Absorbing Boundary Condition (ABC)~\cite{Webb1989}, Boundary Integral (BI)~\cite{Sheng1998}, Perfectly Matched Layer (PML)~\cite{Pekel1995} are required. Zoltan J. Cendes et al.~\cite{Sheng1998} presented the implementation of the tangential vector finite element method to analyze the dielectric waveguides. In the literature various methods like method of moment, spectral-domain methods, finite difference and finite element are adopted for analysis of dielectric waveguide problems. 
\par      J P Webb~\cite{Webb2005} broadly discussed the useful properties of edge elements. Vector finite element or edge element~\cite{Nedelec1980},~\cite{Boffi1999},~\cite{Boffi2001},~\cite{Cendes1991},~\cite{Reddy1994} was proposed to circumvent the spurious solution problem of nodal FEM. Electromagnetic radiation and scattering problems require special elements where normal discontinuity and tangential continuity exists across material interfaces met by edge elements. In~\cite{Barton1987}, M. L. Barton et al. showed that continuity of tangential components of the vector field is sufficient in vector based FEM to compute the magnetic fields. Another advantage of edge elements is that electric or magnetic fields can be directly computed without any differentiation on potentials. Also, no penalty or regularization term is required in edge element framework. While modeling sharp, perfectly conducting objects, the electric field has to be infinite inside the domain and its direction changes rapidly at the sharp edges and corners.  In order to do eigen analysis of such geometries in nodal framework singular trial functions are required, whereas due to tangential continuity of edge elements, no such function is required.
\par      T. V. Yioultsis et.\ al.~\cite{Yioultsis1996}-\cite{Yioultsis1997} presented the systematic approach to construct the higher order tetrahedral edge elements and used them to solve the waveguide problems with material discontinuity. General expressions for the shape functions were presented and unknown coefficients were found by following the decoupling procedure. In~\cite{Rapetti2009}, the formation of higher order Whitney p-elements was shown. In~\cite{Lee1991}, the construction of higher order two dimensional and three dimensional H1 curl elements are presented to solve the electromagnetic scattering problems. In~\cite{Graglia1997}, Roberto D. Graglia et al. presented the general approach of interpolatory vector basis functions of various two dimensional and three dimensional elements.  Edge elements can be constructed by using hierarchial vector basis functions also.  Jon P. Webb~\cite{Webb1999} proposed hierarchical vector basis functions for higher order triangle and tetrahedral finite elements. In~\cite{Seung-CheolLee2003}, Seung-Cheol Lee et al. implemented higher order hierarchical vector finite elements in the field of microwave engineering to the waveguiding structures. In hierarchial type implementation, there can be p refinement in some part of the domain and in some part we can have h refinement. But, in interpolatory type we can only have one type of refinement in the entire domain. To the best knowledge of authors, in all the above literatures, the detailed conversion strategy from nodal FE input file to edge element data-structure is not available. As most of the available mesh generator packages are based on nodal FEM, it will be very useful if such conversion algorithm is developed. In the current work, we have presented a systematic and thorough conversion algorithm.

\par     The remaining article is organized as follows: In section~\ref{edge_data_structure}, we have presented the algorithm in minute details with associated flowcharts and simple conversion examples for different elements like four edge quadrilateral, three edge triangle, twelve edge quadrilateral and eight edge triangle. In section~\ref{numerical_examples}, we have validated our conversion algorithm with several benchmark numerical examples, including all possible complexities like curved surfaces, sharp edges and corners. We have compared our results with available analytical and benchmark solution from literatures.

\section{Conversion algorithm for creating Edge data structure from Nodal data structure}
\label{edge_data_structure}

    For electromagnetic analysis in edge element framework, element data is required in the form of edge data structure. To achieve this, a standalone conversion algorithm is needed where the edge information is generated as output from the supplied nodal data as input. Here, every edge is generated by joining the two nodes of the element. In this section we are presenting such conversion algorithm to different edge elements. Such kind of program is necessary because most of the available commercial mesh generator create and generate nodal element data structure. 

\subsection{Calculation of $\parder{\xi}{x}$, $\parder{\xi}{y}$, $\parder{\eta}{x}$ and $\parder{\eta}{y}$}
 The edge shape functions contains the terms of inverse Jacobian ($\bf{\Gamma}$) and its components $\parder{\xi}{x}$, $\parder{\xi}{y}$, $\parder{\eta}{x}$ and $\parder{\eta}{y}$ which can be obtained as discussed below. From the relation

\begin{equation*}
\begin{Bmatrix}
\parder{f}{\xi} \\[2mm]
\parder{f}{\eta}
\end{Bmatrix} = \bJ \begin{Bmatrix}
\parder{f}{x} \\[2mm]
\parder{f}{y}
\end{Bmatrix} \implies
\begin{Bmatrix}
\parder{f}{x} \\[2mm]
\parder{f}{y}
\end{Bmatrix}=\bJ^{-1}\begin{Bmatrix}
\parder{f}{\xi} \\[2mm]
\parder{f}{\eta}
\end{Bmatrix}  \label{jacobian}
\end{equation*}
 \text{where} $\bJ$ \text{is Jacobian can be written as}
\begin{align*}
	\bJ&=\begin{bmatrix}
J_{11} & J_{12} \\[2mm] J_{21}  & J_{22}
\end{bmatrix}=
\begin{bmatrix}
\parder{x}{\xi} &
\parder{y}{\xi}  \\[2mm]
\parder{x}{\eta} &
\parder{y}{\eta}
\end{bmatrix}
\text{and assume } \bJ^{-1}=\mathbf{\Gamma}=\begin{bmatrix}
\Gamma_{11} & \Gamma_{12} \\[2mm] \Gamma_{21}  & \Gamma_{22}
\end{bmatrix};\\
 \therefore \begin{Bmatrix}
\parder{\xi}{x} \\[2mm]
\parder{\xi}{y}
 \end{Bmatrix}&=
\begin{Bmatrix}
\Gamma_{11} \\[2mm]
\Gamma_{21}
\end{Bmatrix} \text{and }
\begin{Bmatrix}
\parder{\eta}{x} \\[2mm]
\parder{\eta}{y}
\end{Bmatrix}=\begin{Bmatrix}
\Gamma_{12} \\[2mm]
\Gamma_{22}
\end{Bmatrix}
\end{align*}

\subsection{Four edge quadrilateral element}
\begin{figure}[pos=ht]
\centering
\begin{subfigure}{0.4\textwidth}
\centering
\includegraphics[width=0.7\textwidth]{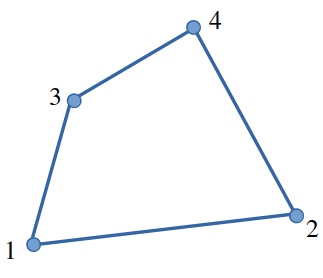}
\caption{Element with four nodes}
\label{loc_nod_con}
\end{subfigure}%
\begin{subfigure}{0.4\columnwidth}
\centering
\includegraphics[width=0.7\textwidth]{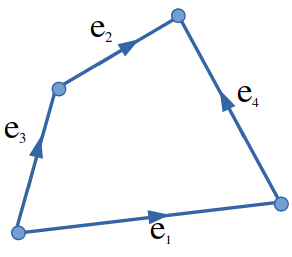}
\caption{Element with four edges}
\label{loc_edge_con}
\end{subfigure}

\caption{Quadrilateral element}
\end{figure}

For four node quadrilateral element Fig.~\ref{loc_nod_con} shows the local nodal connectivity for which Fig.~\ref{loc_edge_con} shows the required local edge connectivity sequence. Edge \textbf{e}$_{1}$ is formed by connecting the local node set \textbf{(1,2)}. Similarly, \textbf{e}$_{2}$, \textbf{e}$_{3}$ and \textbf{e}$_{4}$ are the other three edges generated by connecting the node sets \textbf{(4,3)}, \textbf{(1,4)} and \textbf{(2,3)} respectively. The edge shape functions of four edges are $\bv_1=\frac{l_1}{4}(1-\eta)\del\xi$, $\bv_2=\frac{l_2}{4}(1+\eta)\del\xi$, $\bv_3=\frac{l_3}{4}(1-\xi)\del\eta$ and $\bv_4=\frac{l_4}{4}(1+\xi)\del\eta$ ~\cite{Jinming-finite}. Here $l_1$, $l_2$, $l_3$ and $l_4$ are the lengths of the edges \textbf{e}$_1$, \textbf{e}$_2$, \textbf{e}$_3$, and \textbf{e}$_4$ respectively. 

\subsubsection{Calculation of $\del\bcross\bE$ with edge element}
The components of $\del\bcross\bE$ can be calculated by using the relation,
\begin{equation*}
\del\bcross\bE =
\begin{bmatrix}
\parder{v_{1y}}{x}-\parder{v_{1x}}{y} & \parder{v_{2y}}{x}-\parder{v_{2x}}{y} & \parder{v_{3y}}{x}-\parder{v_{3x}}{y} & \parder{v_{4y}}{x}-\parder{v_{4x}}{y}
\end{bmatrix}
\begin{Bmatrix}
E_1 \\[2mm] E_2 \\[2mm] E_3 \\[2mm] E_4
\end{Bmatrix} =\bB\bE
\end{equation*} where $E_1$, $E_2$, $E_3$ and $E_4$ are the tangential electric fields along the edges 1, 2, 3 and 4 edges respectively. The components of the $\bB$-matrix is obtained as below.
\begin{equation*}
\begin{bmatrix}
\parder{v_{1x}}{x} & \parder{v_{1x}}{y} 
\end{bmatrix} =
\begin{bmatrix}
\parder{v_{1x}}{\xi} & \parder{v_{1x}}{\eta} \\[2mm]
\end{bmatrix}
\begin{bmatrix}
\Gamma_{11} & \Gamma_{21} \\[2mm]
\Gamma_{12} & \Gamma_{22}
\end{bmatrix} 
\end{equation*}  We can derive partial derivative terms $\parder{v_{1x}}{\xi}$, $\parder{v_{1y}}{\eta}$ ... etc with finite difference or by mathematica software~\cite{mathematica}. 
\subsubsection{Conversion algorithm: Generation of edge connectivity array}
\label{generating_edge_connectivity_array}
\begin{figure}[pos=ht]
\centering
\begin{subfigure}{0.4\textwidth}
\centering
\includegraphics[width=0.7\textwidth]{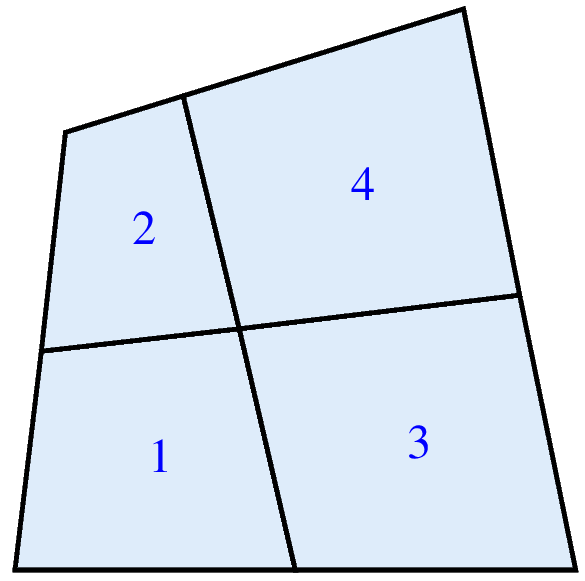}
\caption{Domain with global elements}
\label{glo_ele_quad}
\end{subfigure}%
\begin{subfigure}{0.4\columnwidth}
\centering
\includegraphics[width=0.7\textwidth]{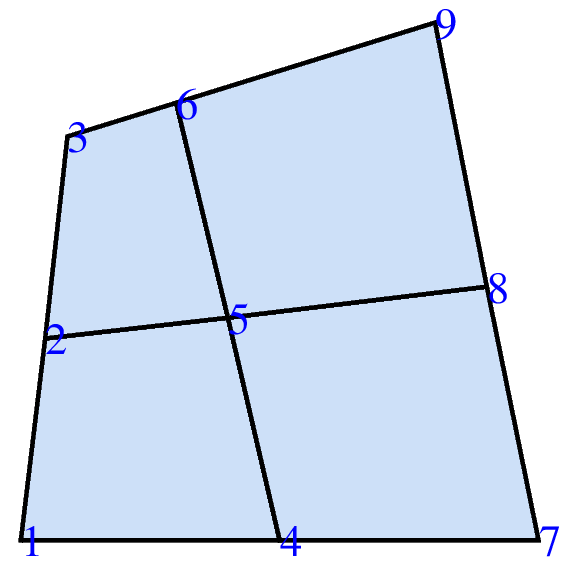}
\caption{Domain with global nodes}
\label{nod_con_quad}
\end{subfigure}
\caption{Discretized domain}
\end{figure}

To describe the conversion algorithm, we have considered a domain which is meshed with quadrilateral elements as shown in Fig.~\ref{glo_ele_quad}. Fig.~\ref{nod_con_quad} shows the global nodal connectivity of the meshed domain. Nodal connectivity list for all the elements for the mesh is given in Table~\ref{elenodcon}. 

\begin{table}[width=.4\linewidth,cols=2,pos=h!]
\caption{Element nodal connectivity} \label{elenodcon}
\begin{tabular*}{\tblwidth}{@{}LL @{} } \toprule
\textbf{Element} &\textbf{Nodal connectivity} \\ 
\textbf{number} & \textbf{(Global node no.)} \\ \midrule
1   &  1, 4, 5, 2 \\
2   &  2, 5, 6, 3 \\
3   &  4, 7, 8, 5 \\
4   &  5, 8, 9, 6 \\ \bottomrule
\end{tabular*}
\end{table}

In this algorithm, the outermost loop runs over the total number of discretized elements and the next inner loop runs over total number of local edges ($nlocedge$) of each element. One global counter $e_l$ is used which is updated to the last assigned global edge number at the end of each element loop. Inside the loop of $nlocedge$, a subroutine \textbf{\emph{edgend}} as shown in Fig.~\ref{edgend} returns global node numbers of two end nodes ($endnd1$ and $endnd2$) of the edge using local nodal connectivity and local edge connectivity information as shown in Fig.~\ref{loc_nod_con} and Fig.~\ref{loc_edge_con}. We follow the convention that each edge direct from $endnd1$ to $endnd2$.

\tikzstyle{startstop} = [rounded rectangle, draw=black, minimum width=2cm, minimum height=0.8cm, text =black, text centered, very thick]
\tikzstyle{io} = [trapezium, trapezium left angle=70, trapezium right angle=110, minimum width=1.3cm, minimum height=1.2cm, text width=6.0cm,text centered, draw=black, text =black,very thick]
\tikzstyle{process} = [rectangle, minimum width=1.3cm, minimum height=1cm, text centered, text width=5.5cm,draw=black, text =black,very thick]
\tikzstyle{decision} = [diamond, aspect=1.6, text centered, text =black,text width=3.5cm,draw=black, very thick]
\tikzstyle{arrow} = [-{Stealth[scale=1.2]},rounded corners,very thick,draw=black,text=black]

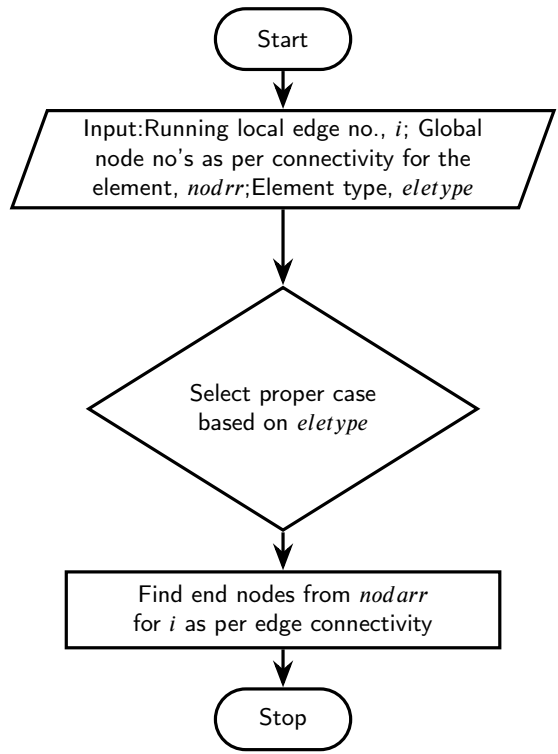
\begin{figure}[pos=h]
\centering
\begin{tikzpicture}[scale=0.4,node distance=0.5cm]
\node (start) [startstop] {Start};
\node (in1) [io, below = of start] {Input:Running local edge no., $i$; Global node no's as per connectivity for the element, $nodrr$;Element type, $eletype$};
\node (dec1) [decision, below = of in1, yshift=-0.5cm] {Select proper case based on $eletype$};
\node (pro1) [process, below = of dec1] {Find end nodes from $nodarr$ for $i$ as per edge connectivity};
\node (stop) [startstop, below = of pro1] {Stop};
\draw [arrow] (start) -- (in1);
\draw [arrow] (in1) -- (dec1);
\draw [arrow] (dec1) -- (pro1);
\draw [arrow] (pro1) -- (stop);
\end{tikzpicture}
\caption{Flow chart of edgend structure}
\label{edgend}
\end{figure}

\begin{table}[width=.8\textwidth,cols=3,pos=h!]
\caption{Local edges and its end nodes of elements of discretized domain.} \label{loc_edge_endnd}
\begin{tabular*}{\tblwidth}{@{}LLL@{}}  \toprule
\textbf{Element} & \textbf{Local edge} & \textbf{End nodes} \\ 
\textbf{number} & \textbf{number} & \textbf{(endnd1, endnd2)} \\ \midrule  
 \multirow{4}{*}{1} & 1 & 1, 4 \\ 
                    & 2 & 2, 5 \\ 
                    & 3 & 1, 2 \\ 
                    & 4 & 4, 5 \\ \midrule 
 \multirow{4}{*}{2} & 1 & 2, 5 \\ 
                    & 2 & 3, 6 \\ 
                    & 3 & 2, 3 \\ 
                    & 4 & 5, 6 \\ \midrule 
 \multirow{4}{*}{3} & 1 & 4, 7 \\  
                    & 2 & 5, 8 \\  
                    & 3 & 4, 5 \\  
                    & 4 & 7, 8 \\ \midrule 
 \multirow{4}{*}{4} & 1 & 5, 8 \\  
                    & 2 & 6, 9 \\  
                    & 3 & 5, 6 \\  
                    & 4 & 8, 9 \\ \bottomrule
\end{tabular*}
\end{table}

	For example, two local nodes 4 and 3 (see Fig.~\ref{loc_nod_con}) are connected to form the local edge \textbf{e}$_{2}$ of the element as shown in Fig.~\ref{loc_edge_con}. These two local node numbers are the position in the nodal connectivity array ($nodarr$) of the element, this array contain global node numbers as shown in Table~\ref{elenodcon}. Thus, for example, for second element, global nodes 6 and 3 are two end nodes for edge \textbf{e}$_{2}$. Local edges and its corresponding connected nodes of the elements for the discretized domain is tabulated in Table~\ref{loc_edge_endnd}.

After successful collection of output from the \textbf{\textit{edgend}} subroutine i.e., information about two end nodes of local edge (\textit{i}), last assigned global edge (\textit{e$_{l}$}) and existing edge connectivity array (\textit{edgearr}) are further supplied into \textbf{\textit{edgedata}} structure as shown in Fig.~\ref{edgedata}. In this data structure at the end of each iteration (for each local edge), different variables like \textit{nodeedgenum, nodeedge, nodeedgexn, edgenode, edgearr} are updated for every endnode (nd1 and nd2) of each edge which are discussed as below.

\tikzstyle{startstop} = [rounded rectangle, draw=black, minimum width=2cm, minimum height=0.8cm, text centered, text=black, very thick]
\tikzstyle{io} = [trapezium, trapezium left angle=70, trapezium right angle=110, minimum width=1.3cm, minimum height=1.2cm, text width=6.0cm,text centered, text=black,draw=black, very thick]
\tikzstyle{process} = [rectangle, minimum width=1.3cm, minimum height=1cm, text centered, text width=5.5cm,draw=black, text=black,very thick]
\tikzstyle{decision} = [diamond, aspect=1.5, text centered, text width=4.5cm,draw=black, text=black,very thick]
\tikzstyle{arrow} = [-{Stealth[scale=1.2]},rounded corners,thick,draw=black,text=black]
\tikzstyle{line} = [-{Stealth[scale=1.2]},thick,draw=black]
\tikzstyle{process1} = [rectangle, minimum width=1.3cm, minimum height=1cm, text centered, text=black,text width=5.5cm,]
\tikzstyle{arr} = [arrow, draw=red]
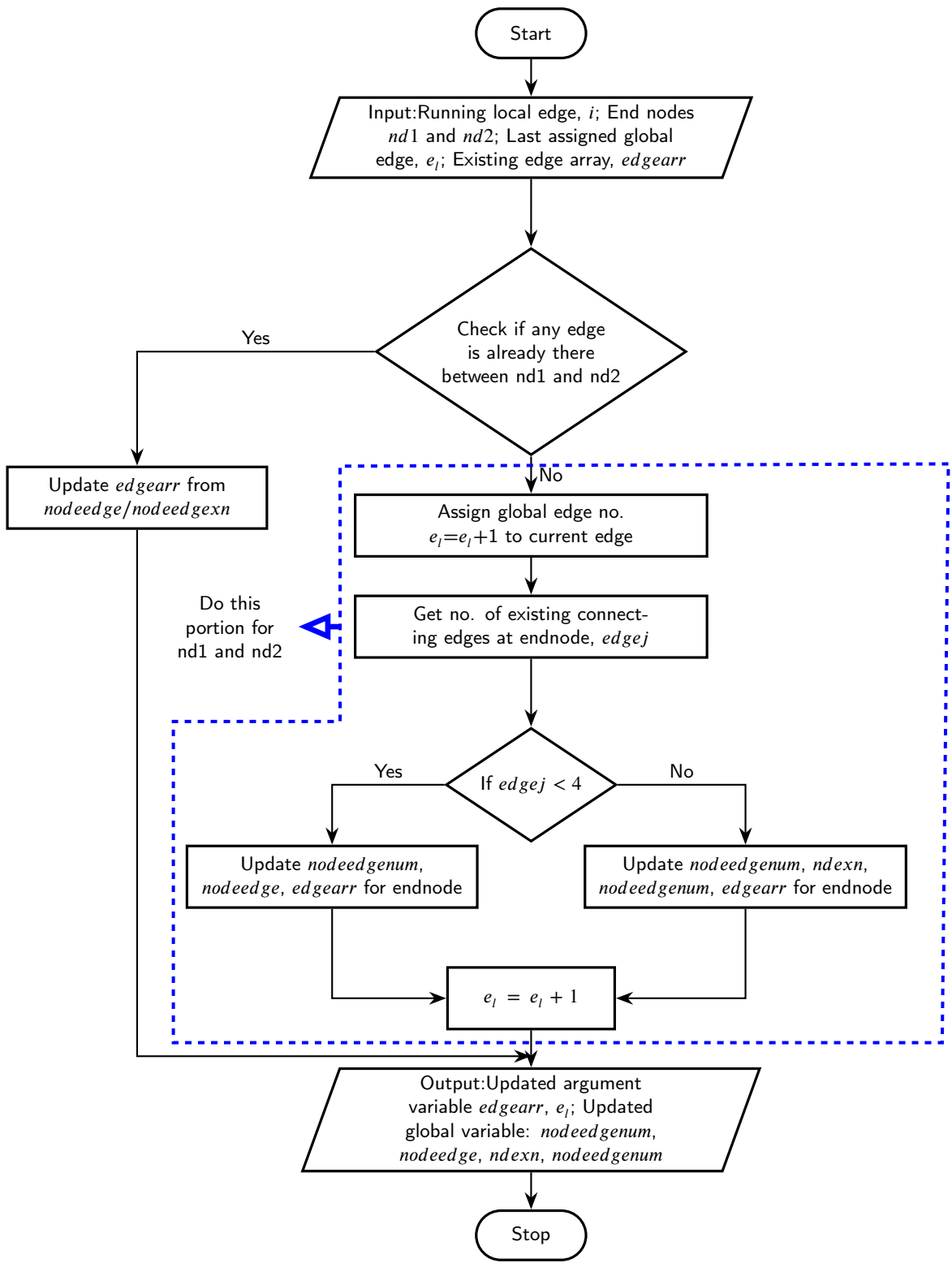
\begin{figure}[pos=h!]
\centering
\begin{tikzpicture}[node distance=0.6cm]
\node (start) [startstop] {Start};
\node (in1) [io, below = of start] {Input:Running local edge, $i$; End nodes $nd1$ and $nd2$; Last assigned global edge, $e_l$; Existing edge array, $edgearr$};
\node (dec1) [decision, below = of in1, yshift=-0.5cm,text width=3cm] {Check if any edge is already there between nd1 and nd2};
\node (pro1) [process, below = of dec1] {Assign global edge no. $e_l$=$e_l$+1 to current edge};
\node (pro2) [process, below = of pro1] {Get no. of existing connecting edges at endnode, $edgej$};
\node (dec2) [decision, below = of pro2, text width = 1.8cm, yshift=-0.5cm] {If $edgej<4$};
\node (pro3) [process,below left = 0.5cm and 0.15cm  of dec2,text width =4.5cm] {Update $nodeedgenum$, $nodeedge$, $edgearr$ for endnode};
\node (pro4) [process,below right = 0.5cm and 0.15cm  of dec2,text width =5cm] {Update $nodeedgenum$, $ndexn$, $nodeedgenum$, $edgearr$ for endnode};
\node (pro5) [process, below = 2cm of dec2,text width=2.5cm, text centered] {$e_l=e_l+1$ };
\node (in2) [io, below = of pro5] {Output:Updated argument variable $edgearr$, $e_l$; Updated global variable: $nodeedgenum$, $nodeedge$, $ndexn$, $nodeedgenum$};
\node (pro6) [process, below left = 1.0cm and 3cm of dec1,text width=4cm] {Update $edgearr$ from $nodeedge/nodeedgexn$};
\node (pro7) [process1, left = 0.9cm of pro2,text width=2cm] {Do this portion for nd1 and nd2};
\node (stop) [startstop, below = of in2] {Stop};
\draw [arrow] (start) -- (in1);
\draw [arrow] (in1) -- (dec1);
\draw [line] (dec1) -- (pro1)node[right,midway,text=black]{No} ;
\draw [arrow] (pro1) -- (pro2);
\draw [arrow] (pro2) -- (dec2);
\draw [line] (dec2) -| (pro3) node [above,pos=0.25,text=black]{Yes};
\draw [line] (dec2) -| (pro4) node [above,pos=0.25,text=black]{No};
\draw [line] (pro4) |- (pro5)  ;
\draw [line] (pro3) |- (pro5)  ;
\draw [line] (dec1) -| (pro6)node[above,pos=0.25,text=black]{Yes} ;
\draw [line] (pro6) |- ($(pro5.south)!0.7!(in2.north)$) ;
\draw [arrow] (pro5) -- (in2);
\draw [arrow] (in2) -- (stop);
\draw [blue, dashed,line width=1.5pt,inner sep=2mm]([shift={( 3.9cm, 0.48cm)}]pro1.north east) -- ([shift={(-2.4cm, 4.7cm)}]dec2.north west) |- ([shift={(-0.2cm, 2.0cm)}]pro3.north west) |-([shift={( 5.3cm,-0.20cm)}]pro5.south east)--cycle;
\draw [-open triangle 45, line width =2.5pt,draw=blue,shorten <=0.2cm](pro2)--(pro7);
\end{tikzpicture}
\caption{Flow chart of edgedata structure}
\label{edgedata}
\end{figure}

\begin{enumerate}
\item $nodeedgenum$:
Global one dimensional static array of dimension maximum number of global nodes (\textit{mx\_nmnode}) in which $i^{th}$ row store the number of edges shared by the $i^{th}$ global node.
\item $nodeedge$:
This is two dimensional array of dimension (\textit{mx\_nmnode} $\bcross$ 8) where 8 columns of $i^{th}$ row store two information of 4 connecting edges shared by $i^{th}$ node. 1\textsuperscript{st} column store global edge no. of 1\textsuperscript{st} connecting edge, second column store other end node of that connecting edge. 3\textsuperscript{rd} and 4\textsuperscript{th} column store (edge no., other end node no.) of the second connecting edge. 5\textsuperscript{th} to 8\textsuperscript{th} column store similar set of information for third and fourth connecting edge. If some node is shared by more than four connecting edges then additional information from 5\textsuperscript{th} edge are stored in \textit{nodeedgeexn}. 
 
   \par Also this array nodeedge have the direction information of the edge. If the $i^{th}$ node is the start node of the edge i.e. if the edge is going from $i^{th}$ node to the other node then edge no. is stored in the odd column as positive integer. If the edge is towards $i^{th}$ node from the other node then edge no. is stored as negative integer. Total no. of negative edges are kept in account with one counter variable.
\item $nodeedgeexn$:
   Global two dimensional array of size(maximum number of nodes in the programme shared by more than 4 edges ($mx\_nmnnd$), 2 $\times$ maximum additional edges after four edges sharing one node ($mx\_exedge$)). Each row consists of the information about the node having more than four connecting edges. It consists of the information from the fifth edge onwards. Global number of the node whose information are stored in the $i^{th}$ row of $nodeedgeexn$ is stored in $i^{th}$ row of $ndexn$. \\
$nndexn$ - No. of nodes which are associated with more than 4 edges. \\
$ndexn(i)$ - $i^{th}$ node number connected with more than 4 edges (according to occurrence) \\
$nodeedgeexn$(i, 1) and (i, 2) are 5\textsuperscript{th} edge no. for the node ndexn(i) and other end node of that edge, (i, 3) and (i, 4) are 6\textsuperscript{th} edge no. for the node ndexn(i) and other end node of that edge and so on.
\item $edgenode$:
Global variable of dimension(Maximum number of edges in the programme, 2) in which $j^{th}$ row contains global number of starting node and ending node of $j^{th}$ edge in two columns respectively.
\item $edgearr$:
Local argument variable in the element loop which store in the process global edge numbers of all the edges of the element according to edge connectivity.
\end{enumerate}

\begin{table*}[width=.49\textwidth,cols=2,pos=ht]
\begin{minipage}{0.49\textwidth}
\caption{Elemental edge connectivity of meshed domain.} \label{edgearr}
\begin{tabular*}{\tblwidth}{@{} LL@{}}  \toprule
\textbf{Element} & \textbf{Edge connectivity}  \\ 
\textbf{number} & \textbf{(Global edge no.)} \\ \midrule  
1 & 1, 2, 3, 4   \\  
2 & 2, 5, 6, 7   \\  
3 & 8, 9, 4, 10  \\  
4 & 9, 11, 7, 12 \\  \bottomrule
\end{tabular*}
\end{minipage}\hfill
\begin{minipage}{0.49\textwidth}
\caption{Edgenode array of nodes of the meshed domain.} \label{edgenodearr}
\begin{tabular*}{\tblwidth}{@{} LLL@{}}  \toprule
\textbf{Global} & \textbf{Starting} & \textbf{End}  \\ 
\textbf{edge(i)} & \textbf{node} & \textbf{node} \\ \midrule  
1 & 1 &  4 \\ 
2 & 2 &  5 \\ 
3 & 1 &  2 \\ 
4 & 4 &  5 \\ 
5 & 3 &  6 \\ 
6 & 2 &  3 \\ 
7 & 5 &  6 \\ 
8 & 4 &  7 \\ 
9 & 5 &  8 \\ 
10 & 7 &  8 \\  
11 & 6 &  9 \\  
12 & 8 &  9 \\ \bottomrule 
\end{tabular*}
\end{minipage}
\end{table*}

\begin{figure}[pos=h]
\centering
\includegraphics[width = 0.33\textwidth]{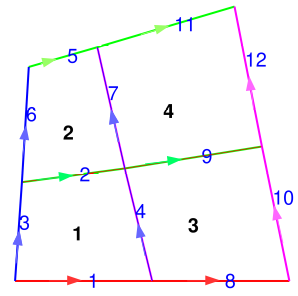}
\caption{Element to edge connectivity sequence of the meshed domain}
\label{fig:edgecon}
\end{figure}

 In \textbf{\textit{edgedata}} subroutine, there is a running counter called $edgej$ which stores the associated number of edges of the starting end nodes of the current local edge. The information is interchanged between local $edgej$ and global array $nodeedgenum$. In order to understand the update of the variables let us first summarize stored values in different variables after completion of element loop for first element.
\begin{enumerate}[(i)]

\item first, second, fourth and fifth rows of `$nodeedgenum$' array are assigned `two' because first element has global nodes 1, 2, 4, and 5 (see Fig.~\ref{nod_con_quad}), and after the loop for element 1 each of these nodes associate with `two' edges.
\item Last assigned global edge i.e., $e_l$ with the value 4.
\item $nodeedge$ will have following values as shown in Table~\ref{nodeedgearr1}.

\begin{table}[width=.8\textwidth,cols=9,pos=h]
\caption{Nodeedge array of nodes after element loop of the first element} \label{nodeedgearr1}
\begin{tabular*}{\tblwidth}{@{} LLLLLLLLL@{}}  \toprule
 &  \multicolumn{8}{c}{\textbf{Associated edge and}} \\
\textbf{Global} & \multicolumn{8}{c}{\textbf{second node of the edge(nodeedge(1:8))}} \\ \cline{2-9}
\textbf{node} & \textbf{1$^{st}$} & \textbf{Other} &  \textbf{2$^{nd}$} & \textbf{Other} & \textbf{3$^{rd}$} & \textbf{Other}  & \textbf{4$^{th}$} & \textbf{Other} \\ 
 & \textbf{edge} & \textbf{node} & \textbf{edge} & \textbf{node} & \textbf{edge} & \textbf{node} & \textbf{edge} & \textbf{node} \\\midrule
1 &   1  & 4 &  3  & 2 & 0  & 0 & 0 & 0 \\ 
2 &  -3  & 1 &  2  & 5 & 0  & 0 & 0 & 0 \\ 
3 &   0  & 0 &  0  & 0 & 0  & 0 & 0 & 0 \\ 
4 &  -1  & 1 &  4  & 5 & 0  & 0 & 0 & 0 \\ 
5 &  -2  & 2 & -4  & 4 & 0  & 0 & 0 & 0 \\ 
6 &   0  & 0 &  0  & 0 & 0  & 0 & 0 & 0 \\ 
7 &   0  & 0 &  0  & 0 & 0  & 0 & 0 & 0 \\ 
8 &   0  & 0 &  0  & 0 & 0  & 0 & 0 & 0 \\ 
9 &   0  & 0 &  0  & 0 & 0  & 0 & 0 & 0 \\ \bottomrule
\end{tabular*}
\end{table}

For first row (for global node 1) we have  1 (associated first edge), 4 (other node of the first edge), 3 (associated second edge), and 2 (other node of the second edge). The reason for the negative sign in (2,1) position of $nodeedge$ is that the first associated edge of the node 2 i.e. edge 3 is directing from the other node (1, stored in (2,2)) towards the current node 2 (see Fig.~\ref{fig:edgecon}). For node 1, both the associated edges (1 and 3) are directing from the current node 1 to the respective other nodes. Therefore, those edge numbers are stored with `+' sign. In the first row (i.e., for first global edge) of $edgenode$ array, starting node (1) and end node (4) are stored as shown in Fig.~\ref{edgenodearr}. Similarly for the other global edges of the element 1, starting and end nodes are stored in 2\textsuperscript{nd}, 3\textsuperscript{rd} and 4\textsuperscript{th} rows of this array. Table~\ref{edgenodearr} shows such information for the edges 1 to 4.
\item First row of edge connectivity array is stored with the edge numbers of the first element, 1, 2, 3 and 4 as shown in the first row of Table~\ref{edgearr}.
\end{enumerate}

   Now for the next entity of the outer element loop i.e., for element 2, for local edge 1 we have starting and end nodes as 2 and 5 (see Fig.~\ref{nod_con_quad}) as per convention of Fig.~\ref{loc_nod_con} and Fig.~\ref{loc_edge_con}. At first, from $nodeedgenum$, we get the total number of already associated edges of the starting node. For our starting node 2, there are two associated edges. After that we get the other node numbers of the associated edges from the even columns of row 2 (our current node) of $nodeedge$ array. If any of these other node number match with our current end node (5) then the associated edge no. (available in respective odd column) will be the global edge no. of that local edge. In this case the edge no. is 2, we will update the current edge number with global edge number 2. In this local edge loop we will not update last assigned global edge($e_l$). Also, the first column of $edgearr$ is assigned with 2.

    \par For the 2\textsuperscript{nd} local edge starting and end nodes are 3 and 6 as shown in Fig.~\ref{nod_con_quad}. As there is no data available in 3\textsuperscript{rd} row of $nodeedge$ array, the end node 3 is appearing for the first time. So, a new edge number is assigned just by updating $e_l$ to $e_l+1$ i.e., with the digit 5. Therefore, in the $nodeedge$ array, first and second columns of third(current end node) row are assigned with 5 (associated edge no.) and 6 (other end node) respectively. For end node number 6, the sixth row of this array is updated with -5 and 3 \text{in} the first two columns. Because the edge 5 is pointing away from the current end node 6, a negative symbol is assigned to the digit 5. 5\textsuperscript{th} row of the $edgenode$ array is updated with the end node informations (3 and 6) of this new edge (5). In the 2\textsuperscript{nd} column of $edgearr$, this new edge number 5 is assigned. In $nodeedgenum$, existing number in third and sixth rows are incremented by 1 because global node 3 and global node 6 become associated with new edge `5' in this local edge loop.

   \par For the third local edge, 2 and 3 are supplied as starting and end nodes (see Fig.~\ref{nod_con_quad}) respectively from the $edgend$ subroutine. For starting node 2, there are two connected edges 3 and 2. Other end nodes of these connected edges (available in the even columns) i.e. 1 and 5 are not matching with other end node 3. As  there is no already existing edge between these two nodes, last assigned global edge ($e_l$) is incremented from 5 to 6. After this, $edgearr$ is updated with this value in the third row. In the second row of $nodeedge$ array this new edge data (global edge 6) and its other end node (3) are updated in the fifth and sixth columns. For end node 3, the corresponding row of this array (3\textsuperscript{rd} row) is updated with -6 and 2 \text{in} the third and fourth columns. `-' symbol is assigned to the digit 6 because the edge 6 is pointing away from the current node 3. Also, the second row of $nodeedgenum$ array(related to 2\textsuperscript{nd} global node) is updated from the previous count 2 to 3. Similarly, third row (related to 3\textsuperscript{rd} global node) is updated from 1 to 2. $edgenode$ array is also updated with the end nodes (2 and 3) of the newly formed edge in the first and second columns of the corresponding (sixth) row.
   \par For the fourth local edge, end nodes as 5 (starting node) and 6 (end node), $e_l$ with 6 and existing $edgearr$ (2, 5 and 6) are supplied. After checking the end node 6 with the other nodes (2 and 4) available in the fifth row of the $nodedge$ array, $e_l$ is updated from 6 to 7. This value is assigned to the running local edge. Now, fourth column of $edgearr$ is updated with 7. Thus we complete the 2\textsuperscript{nd} row of Table~\ref{edgearr} which shows the element to edge connectivity array of second element. Fig.~\ref{fig:edgecon} shows such element to edge connectivity for all the global elements of finite element meshed domain.
   
\tikzstyle{startstop} = [rounded rectangle, draw=black, minimum width=2cm, minimum height=0.8cm, text centered, text=black,very thick]
\tikzstyle{io} = [trapezium, trapezium left angle=70, trapezium right angle=110, minimum width=1.3cm, minimum height=1.2cm, text width=6.0cm,text centered, text=black,draw=black, very thick]
\tikzstyle{process} = [rectangle, minimum width=1.3cm, minimum height=0.5cm, text centered, text width=5.5cm,draw=black, text=black,very thick]
\tikzstyle{decision} = [diamond, aspect=1.5, text centered, text width=4.5cm,draw=black, text=black,very thick]
\tikzstyle{arrow} = [-{Stealth[scale=1.2]},rounded corners,thick,draw=black,text=black]
\tikzstyle{line} = [-{Stealth[scale=1.2]},thick,draw=black]
\tikzstyle{process1} = [rectangle, minimum width=1.3cm, minimum height=1cm, text centered, text=black,text width=5.5cm,]
\tikzstyle{arr} = [arrow, draw=red]

\begin{figure}[pos=h!]
\centering
\begin{tikzpicture}[node distance=0.6cm]
\node (start) [startstop] {Start};
\node (in1) [io, below = of start] {Input:Nodal connectivity array, $nodecon$; Element type, $eletype$; No. of elements, $nele$};
\node (pro1) [process, below = of in1,text width = 1.8 cm] {ele = 1};
\node (pro2) [process, below = of pro1] {Based on $eletype$, find $nodarr$ from $nodecon$; No. of edges in each element, $nedge$};
\node (pro3) [process, below = of pro2,text width = 1.8cm ] {i = 1};
\node (pro4) [process, below  = of pro3,text width =6cm] {Find $nd_1$ and $nd_2$ using \textbf{\textit{edgend}}};
\node (pro5) [process, below  = of pro4,text width =6cm] {Update $edgearr(el)$, $nodeedgenum$, $nodeedge$, $ndexn$, $nodeedgeexn$ using \textbf{\textit{edgedata}}};
\node (pro6) [process, below =  of pro5,text width=2cm] {i = i + 1 };
\node (pro7) [process, below =  of pro6,text width=6cm,text centered] {Update edge connectivity array, $edgecon$ appending $edgearr$ to it};
\node (pro8) [process, below =  of pro7,text width=3cm] {ele = ele + 1};
\node (in2) [io, below = of pro8] {Output:Updated $edgecon$, $nodeedgenum$, $nodeedge$, $ndexn$, $nodeedgeexn$};
\node (stop) [startstop, below = of in2] {Stop};

\draw [arrow] (start) -- (in1);
\draw [arrow] (in1) -- (pro1);
\draw [arrow] (pro1) -- (pro2);
\draw [arrow] (pro2) -- (pro3);
\draw [arrow] (pro3) -- (pro4);
\draw [arrow] (pro4) -- (pro5);
\draw [arrow] (pro5) -- (pro6);
\draw [arrow] (pro6) -- (pro7);
\draw [line] (pro6) -| ++(-3.5cm,0.5cm) |-node[sloped, above,pos=0.25,text=black]{Loop i = 1 to nedge}  ($(pro3.south)!0.4!(pro4.north)$);
\draw [line] (pro8) -| ++(-5.5cm,1.0cm) |-node[sloped, above,pos=.25,text=black]{Loop ele = 1 to nele} ($(pro1.south)!0.4!(pro2.north)$);
\draw [arrow] (pro7) -- (pro8);
\draw [arrow] (pro8) -- (in2);
\draw [arrow] (in2) -- (stop);
\end{tikzpicture}
\caption{Flow chart of node to edge structure}
\label{node2edge}
\end{figure}
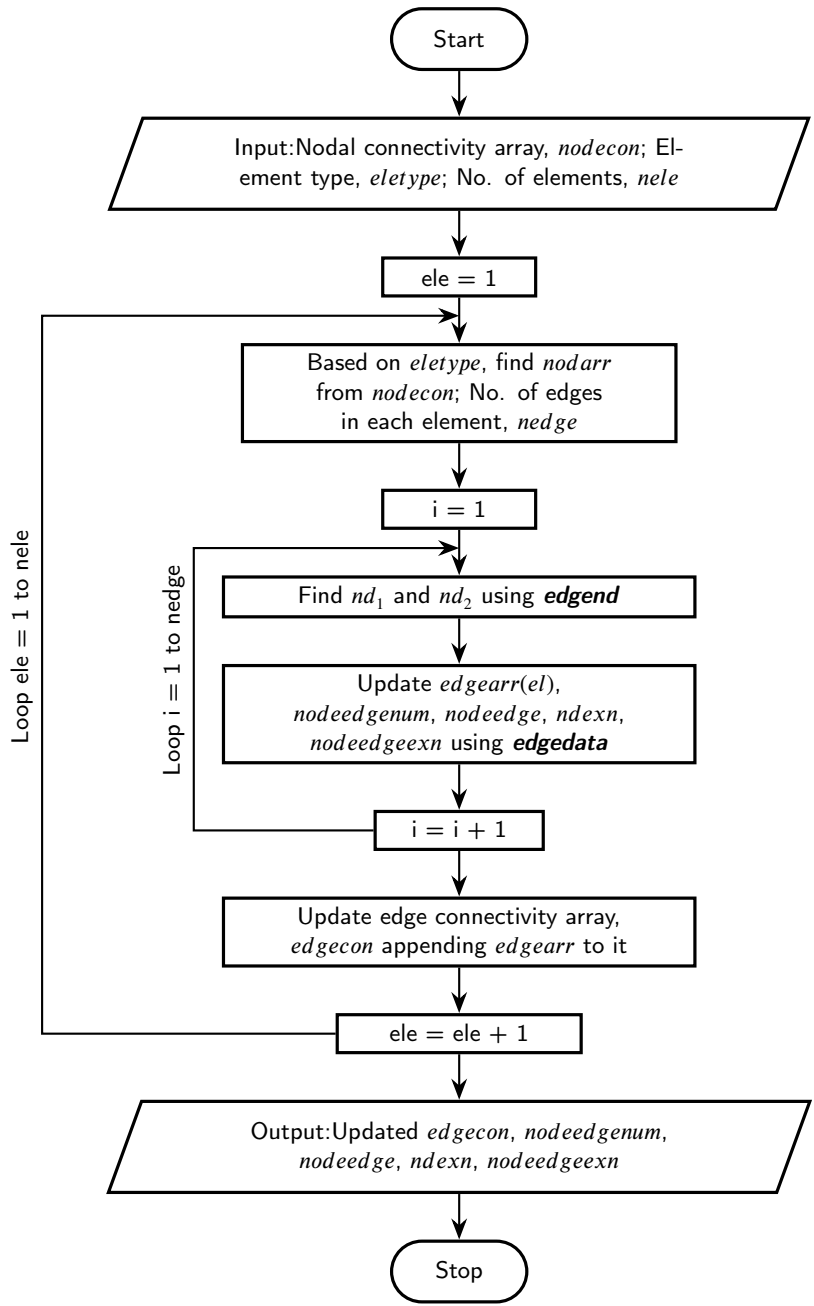
   
   After this, $nodeedge$ array is updated with the new edge number 7 and its associated other node 6 \text{in} the 5\textsuperscript{th} and 6\textsuperscript{th} columns of fifth row. Similarly, in the 6\textsuperscript{th} row corresponding to global node 6, third and fourth columns is updated with the digits -7 and 5 respectively. Also, the existing no. of edges in the fifth and sixth rows of $nodeedgenum$ array are incremented by 1 and updated as 3 and 2. Finally, global edge number 7 and its end nodes (5 and 6) are stored corresponding row (seventh) of $edgenode$ array as shown in Table~\ref{edgenodearr}.

\begin{table}[width=.99\textwidth,cols=10,pos=ht]
\caption{Nodeedge array of the global nodes of the meshed domain.} \label{nodeedgearr}
\begin{tabular*}{\tblwidth}{@{} LLLLLLLLLL@{}}  
\toprule
 & \textbf{Total no. of} & \multicolumn{8}{c}{\textbf{Associated edge and}} \\ 
\textbf{Global} & \textbf{connecting edges} & \multicolumn{8}{c}{\textbf{second node of the edge(nodeedge(1:8))}} \\ \cline{3-10}  
\textbf{node} & \textbf{(nodeedgenum)} & \textbf{1$^{st}$} & \textbf{Other} &  \textbf{2$^{nd}$ } & \textbf{Other} & \textbf{3$^{rd}$} & \textbf{Other}  & \textbf{4$^{th}$} & \textbf{Other} \\ 
 & \textbf{} & \textbf{edge} & \textbf{node} & \textbf{edge} & \textbf{node} &  \textbf{edge} & \textbf{node} & \textbf{edge} & \textbf{node} \\\midrule
1 & 2 &  1  & 4 &  3  & 2 & -  & - & - & - \\  
2 & 3 & -3  & 1 &  2  & 5 & 6  & 3 & - & - \\  
3 & 2 &  5  & 6 & -6  & 2 & -  & - & - & - \\  
4 & 3 & -1  & 1 &  4  & 5 & 8  & 7 & - & - \\  
5 & 4 & -2  & 2 & -4  & 4 & 7  & 6 & 9 & 8 \\  
6 & 3 & -5  & 3 & -7  & 5 & 11 & 9 & - & - \\  
7 & 2 & -8  & 4 & 10  & 8 & -  & - & - & - \\  
8 & 3 & -9  & 5 & -10 & 8 & 12 & 9 & - & - \\  
9 & 2 & -11 & 6 & -12 & 8 & -  & - & - & - \\ \bottomrule 
\end{tabular*}
\end{table}

 After updation of all the global variables, program comes out from the inner (local edge) loop as shown in Fig.~\ref{node2edge} and enters in to outer (element) loop after incrementing as $ele=ele+1$. Now, the program runs for local edges of third element. This process repeats until all the discretized elements are finished. Table~\ref{edgearr} shows the edge connectivity array for the entire domain; Table~\ref{edgenodearr} shows the end nodes for all the edges; Table~\ref{nodeedgearr} shows complete $nodeedge$ array of all the global nodes of the finite element meshed domain after finishing the outer loop for all four elements.

\par This algorithm can be implemented to other edge elements. These edge elements include 3-edge triangle, 8-edge triangle and 12-edge quadrilateral elements. In the following sections we discuss about the implementation of these elements.

\subsection{Three edge triangular element}

\begin{figure}[pos=h!]
\centering
\begin{subfigure}{0.4\textwidth}
\centering
\includegraphics[width=0.7\textwidth]{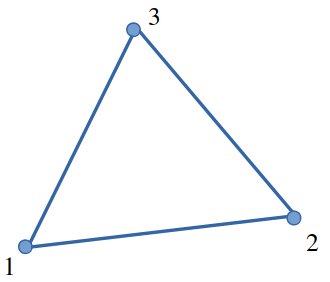}
\caption{Element with three nodes}
\label{3_nod_tri}
\end{subfigure}%
\begin{subfigure}{0.4\columnwidth}
\centering
\includegraphics[width=0.7\textwidth]{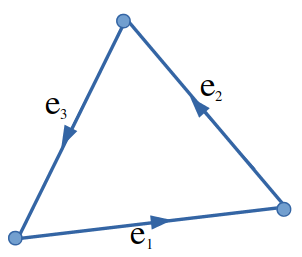}
\caption{Element with three edges}
\label{3_edge_tri}
\end{subfigure}
\caption{Triangular element}
\end{figure}

\begin{figure}[pos=h!]
\centering
\begin{subfigure}{0.5\textwidth}
\centering
\includegraphics[width=0.55\textwidth]{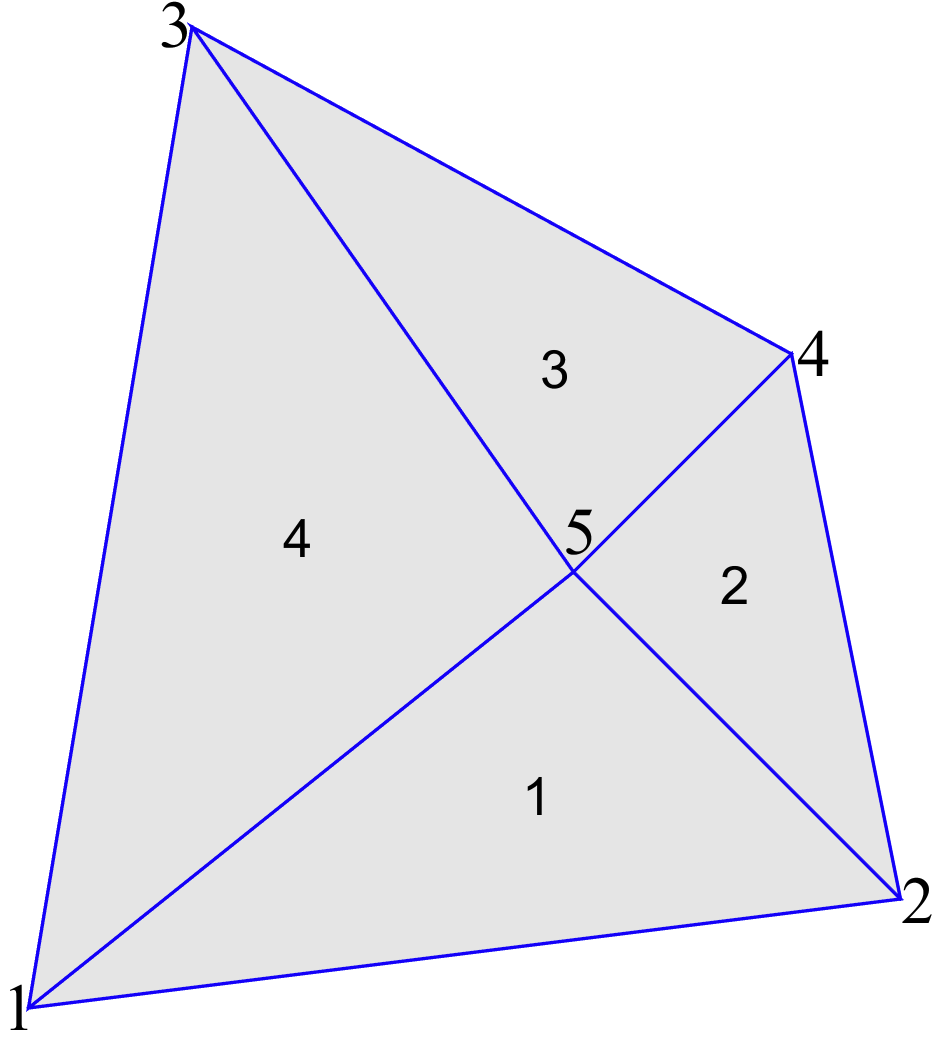}
\caption{Domain discretized with 3 node triangular element}
\label{quadrilateral_mesh_t3_nodes}
\end{subfigure}%
\begin{subfigure}{0.5\columnwidth}
\centering
\includegraphics[width=0.5\textwidth]{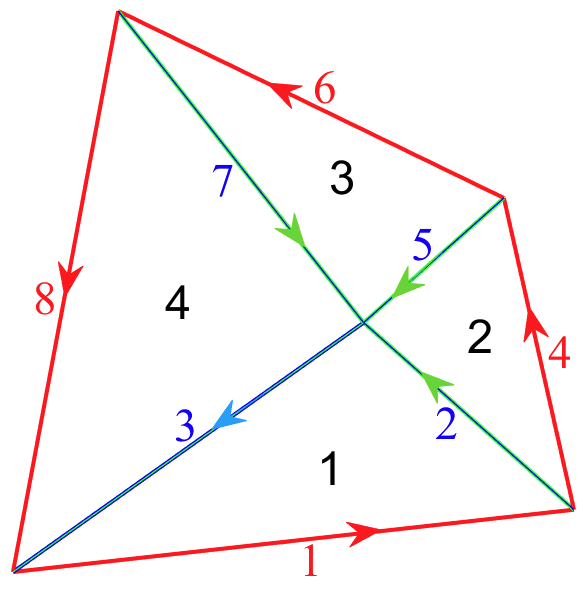}
\caption{Domain discretized with 3 edge triangular element}
\label{quadrilateral_mesh_t3_edges}
\end{subfigure}
\caption{Discretized domain}
\end{figure}

Lower order three edge triangular element is formed from three node triangle and Fig.~\ref{3_nod_tri} shows the nodal connectivity (local) of the triangular element. For this element the expected edge connectivity can be shown in Fig.~\ref{3_edge_tri}. With the help of node sets \textbf{(1,2)}, \textbf{(2,3)} and \textbf{(3,1)} three edges \textbf{e}$_{1}$, \textbf{e}$_{2}$ and \textbf{e}$_{3}$ are formed respectively. Fig.~\ref{quadrilateral_mesh_t3_nodes} shows one general domain discretized with nodal elements which can be transformed into a domain discretized with three edge triangular elements as shown in Fig.~\ref{quadrilateral_mesh_t3_edges}. For this triangular element three edge shape functions are~\cite{Jinming-finite} $\bv_1=l_1(\xi\del \eta-\eta\del \xi)$, ${\bv}_2=l_2(-\eta\del \xi-(1-\xi)\del \eta)$ and $\bv_3=l_3((1-\eta)\del \xi+\xi\del \eta)$ where $l_1$, $l_2$ and $l_3$ are the edge lengths of three edges. 

\subsubsection{Calculation of $\del\bcross\bE$}
From the relation,
\begin{align*}
\del\bcross\bE & =
\begin{bmatrix}
\parder{v_{1y}}{x}-\parder{v_{1x}}{y} & \parder{v_{2y}}{x}-\parder{v_{2x}}{y} & \parder{v_{3y}}{x}-\parder{v_{3x}}{y}
\end{bmatrix}
\begin{Bmatrix}
E_1 \\ E_2 \\ E_3
\end{Bmatrix}
= \bB\bE
\end{align*}
where $E_1$, $E_2$ and $E_3$ are tangential components of electric fields along the three edges \textbf{e}$_1$, \textbf{e}$_2$ and \textbf{e}$_3$ respectively.
The components of the $\bB$-matrix can be obtained by using the relation,
\begin{equation*}
\begin{bmatrix}
\parder{v_{1x}}{x} & \parder{v_{1x}}{y} \\[2mm]
\end{bmatrix} =
\begin{bmatrix}
\parder{v_{1x}}{\xi} & \parder{v_{1x}}{\eta} \\[2mm]
\end{bmatrix}
\begin{bmatrix}
\Gamma_{11} & \Gamma_{21} \\[2mm]
\Gamma_{12} & \Gamma_{22}
\end{bmatrix} 
\end{equation*}
 We can derive the $\parder{v_{1x}}{\xi}$, $\parder{v_{1x}}{\eta}$ ... etc explicitly.

\subsubsection{Conversion algorithm}
\label{nodeedgeexn_circle_domain}

\begin{figure}[pos=h!]
\centering
\begin{subfigure}{0.5\textwidth}
\centering
\includegraphics[width = 0.6\textwidth]{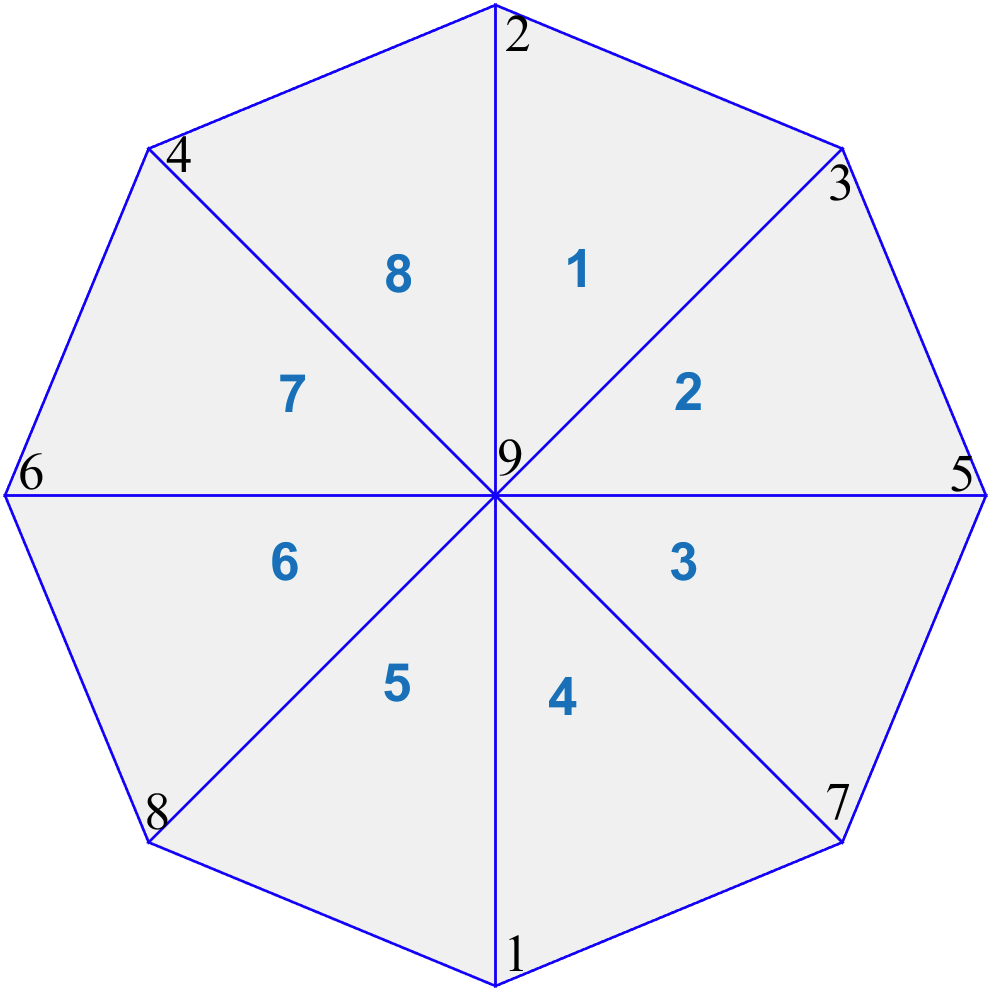}
\caption{Nodal connectivity of the elements}
\label{node_coarse_cir_dom}
\end{subfigure}%
\begin{subfigure}{0.5\textwidth}
\centering
\includegraphics[width = 0.6\textwidth]{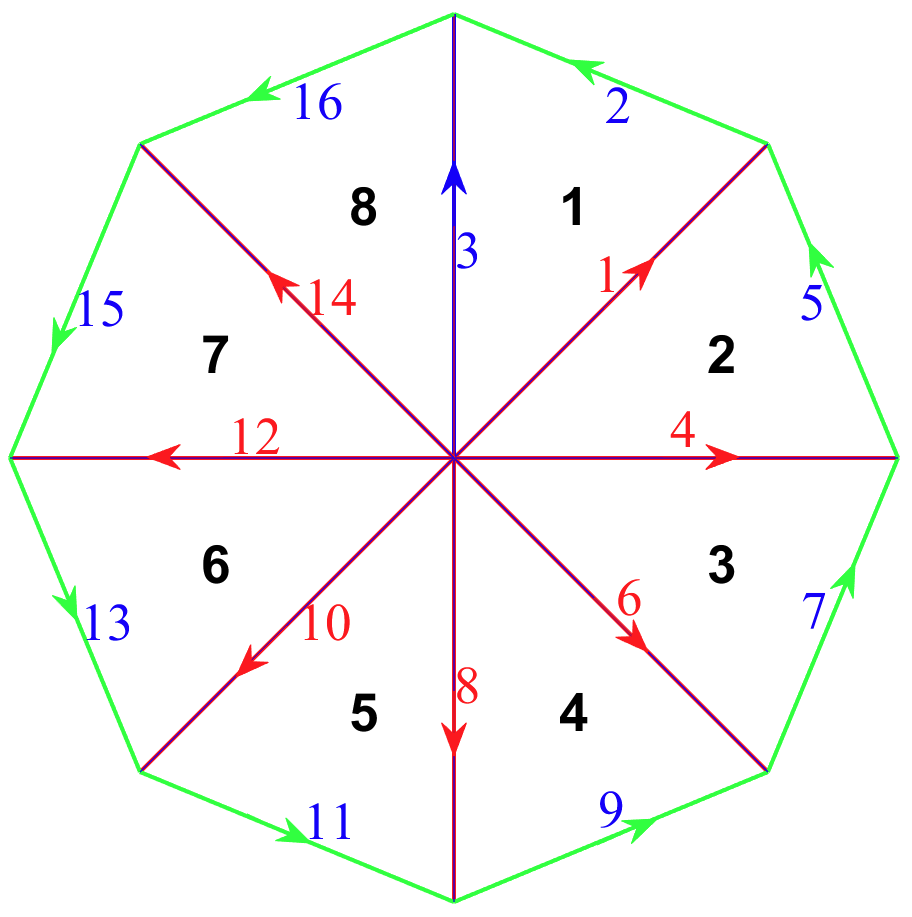}
\caption{Edge connectivity of the elements}
\label{edge_coarse_cir_dom}
\end{subfigure}
\caption{Discretized circular domain with coarse mesh}
\end{figure}

Like four edge quadrilateral elements, edge data structure for three edge triangle can be constructed as described in section \ref{generating_edge_connectivity_array}. If any node is shared by more than four edges then some additional data-structure is required as follows. In this algorithm if any $i^{th}$ node is shared by more than 4 edges then the information from the fifth edge onwards is stored in $nodeedgeexn$, $nndexn$ and $ndexn$ arrays. In entire geometry, there will be few nodes which will be associated with more than four edges. Our \textit{nodeedge} array has number of rows as total number of nodes and number of column as 8 \text{in} order to accommodate first four associated edges which is very common. \textit{nodeedgeexn} is initialized with number of rows far less than total number of nodes, it has 16 columns to accommodate next eight edges. When we come across a node which is associated with more than four edges \textit{nndexn} is incremented by 1. Suppose in a flow, we have jth occurrence of such node which is associated with more than four edges. Then \textit{ndexn(j)} will store the corresponding node number. $j^{th}$ row of \textit{nodeedgeexn} will store the information of associated edge and other node from the fifth edge onwards. This can be understood from the case as shown in Fig.~\ref{node_coarse_cir_dom} and Fig.~\ref{edge_coarse_cir_dom} Node no. 9 of the discretized domain is shared by 8 edges. So the information of the first four connecting edges of Node no. 9 and its other end nodes are stored in $nodeedge$ array as shown in Table~\ref{nodeedgearr_cir_dom} after element loop for fifth element. But from the fifth to eighth associated edge of node 9, the information is stored in $nodeedgeexn$ array. Table~\ref{nodeedgeexnarr5} shows updation of the information of $nodeedgeexn$ array after the end of the element loop for fifth element. After the element loop runs for the remaining existing elements, $nodeedgeexn$ array is updated as shown in Table~\ref{nodeedgeexnarr}. Here, \textit{nndexn} will be 1 and $ndexn(1)$ will be 9. We have solved a numerical example of similar domain as discussed in section~\ref{circular_domain} and section~\ref{cracked_circular_domain} to obtain eigenvalues.

\begin{table}[width=.99\textwidth,cols=10,pos=ht]
\caption{Nodeedge array of nodes after element loop of the fifth element.} \label{nodeedgearr_cir_dom}
\begin{tabular*}{\tblwidth}{@{} LLLLLLLLLL@{}}  
\toprule
 & \textbf{Total no. of} & \multicolumn{8}{c}{\textbf{Associated edge and}} \\ 
\textbf{Global} & \textbf{connecting edges} & \multicolumn{8}{c}{\textbf{second node of the edge(nodeedge(1:8))}} \\ \cline{3-10}  
\textbf{node} & \textbf{(nodeedgenum)} & \textbf{1$^{st}$} & \textbf{Other} &  \textbf{2$^{nd}$} & \textbf{Other} & \textbf{3$^{rd}$} & \textbf{Other}  & \textbf{4$^{th}$} & \textbf{Other} \\ 
 & \textbf{} & \textbf{edge} & \textbf{node} & \textbf{edge} & \textbf{node} &  \textbf{edge} & \textbf{node} & \textbf{edge} & \textbf{node} \\\midrule
9 & 6  &  1  & 3 &  3  & 2 & 4  & 5 & 6 & 7 \\ 
2 & 2  & -2  & 3 & -3  & 9 & -  & - & - & - \\ 
3 & 3  & -1  & 9 &  2  & 2 &-5  & 5 & - & - \\ 
5 & 3  & -4  & 9 &  5  & 3 &-7  & 7 & - & - \\ 
7 & 3  & -6  & 9 &  7  & 5 &-9  & 1 & - & - \\ 
1 & 3  & -8  & 9 &  9  & 7 &-11 & 8 & - & - \\ 
8 & 2  & -10 & 9 & 11  & 1 & -  & - & - & - \\ 
- & -  &  -  & - &  -  & - & -  & - & - & - \\ \bottomrule
\end{tabular*}
\end{table}

\begin{table}[width=.99\textwidth,cols=17,pos=ht]
\caption{Nodeedgeexn array of nodes after element loop of the fifth element} \label{nodeedgeexnarr5}
\begin{tabular*}{\tblwidth}{@{} LLLLLLLLLLLLLLLLL@{}}  \toprule
 &  \multicolumn{16}{c}{\textbf{Associated edge and }} \\ 
\textbf{Global}& \multicolumn{16}{c}{\textbf{second node of the edge (nodeedgeexn(1:16))}} \\ \cline{2-17}
\textbf{node} & \textbf{5$^{th}$} & {\textbf{other}} & \textbf{6$^{th}$} & \textbf{other} & \textbf{7$^{th}$} & \textbf{other} & \textbf{8$^{th}$} & \textbf{other} & \textbf{9$^{th}$} & \textbf{other} & \textbf{10$^{th}$} & \textbf{other} & \textbf{11$^{th}$} & \textbf{other} & \textbf{12$^{th}$} & \textbf{other}  \\ 
 & \textbf{edge} & \textbf{node} & \textbf{edge} & \textbf{node} & \textbf{edge} & \textbf{node} & \textbf{edge} & \textbf{node} & \textbf{edge} & \textbf{node} & \textbf{edge} & \textbf{node} & \textbf{edge} & \textbf{node} & \textbf{edge} & \textbf{node}  \\ \hline
9 &   8  & 1 &  10  & 8 & -  & - & - & - & - & - & - & - & - & - & - & - \\ 
- &   -  & - &  -   & - & -  & - & - & - & - & - & - & - & - & - & - & - \\ \bottomrule
\end{tabular*}
\end{table}


\begin{table}[width=.99\textwidth,cols=17,pos=ht]
\caption{Nodeedgeexn array of nodes of the discretized domain after complete conversion.} \label{nodeedgeexnarr}
\begin{tabular*}{\tblwidth}{@{} LLLLLLLLLLLLLLLLL@{}}  \toprule
 &  \multicolumn{16}{c}{\textbf{Associated edge and }} \\ 
\textbf{Global}& \multicolumn{16}{c}{\textbf{second node of the edge (nodeedgeexn(1:16))}} \\ \cline{2-17}
\textbf{node} & \textbf{5$^{th}$} & {\textbf{other}} & \textbf{6$^{th}$} & \textbf{other} & \textbf{7$^{th}$} & \textbf{other} & \textbf{8$^{th}$} & \textbf{other} & \textbf{9$^{th}$} & \textbf{other} & \textbf{10$^{th}$} & \textbf{other} & \textbf{11$^{th}$} & \textbf{other} & \textbf{12$^{th}$} & \textbf{other}  \\ 
 & \textbf{edge} & \textbf{node} & \textbf{edge} & \textbf{node} & \textbf{edge} & \textbf{node} & \textbf{edge} & \textbf{node} & \textbf{edge} & \textbf{node} & \textbf{edge} & \textbf{node} & \textbf{edge} & \textbf{node} & \textbf{edge} & \textbf{node}  \\ \hline
9 &   8  & 1 &  10  & 8 & 12  & 6 & 14 & 4 & - & - & - & - & - & - & - & - \\ 
- &   -  & - &  -   & - & -  & - & - & - & - & - & - & - & - & - & - & - \\ \bottomrule
\end{tabular*}
\end{table}

\subsection{Eight edge triangular element}

\begin{figure}[pos=h!]
\centering
\begin{subfigure}{0.5\textwidth}
\centering
\includegraphics[width=0.6\textwidth]{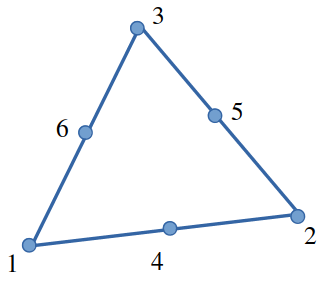}
\caption{Triangular element with six nodes}
\label{6_nod_tri}
\end{subfigure}%
\begin{subfigure}{0.5\columnwidth}
\centering
\includegraphics[width=0.6\textwidth]{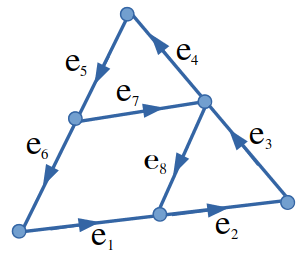}
\caption{Triangular element with eight edges}
\label{8_edge_tri}
\end{subfigure}
\caption{Higher order triangular element}
\end{figure}
\begin{figure}[pos=h!]
\centering
\begin{subfigure}{0.5\textwidth}
\centering
\includegraphics[width=0.5\textwidth]{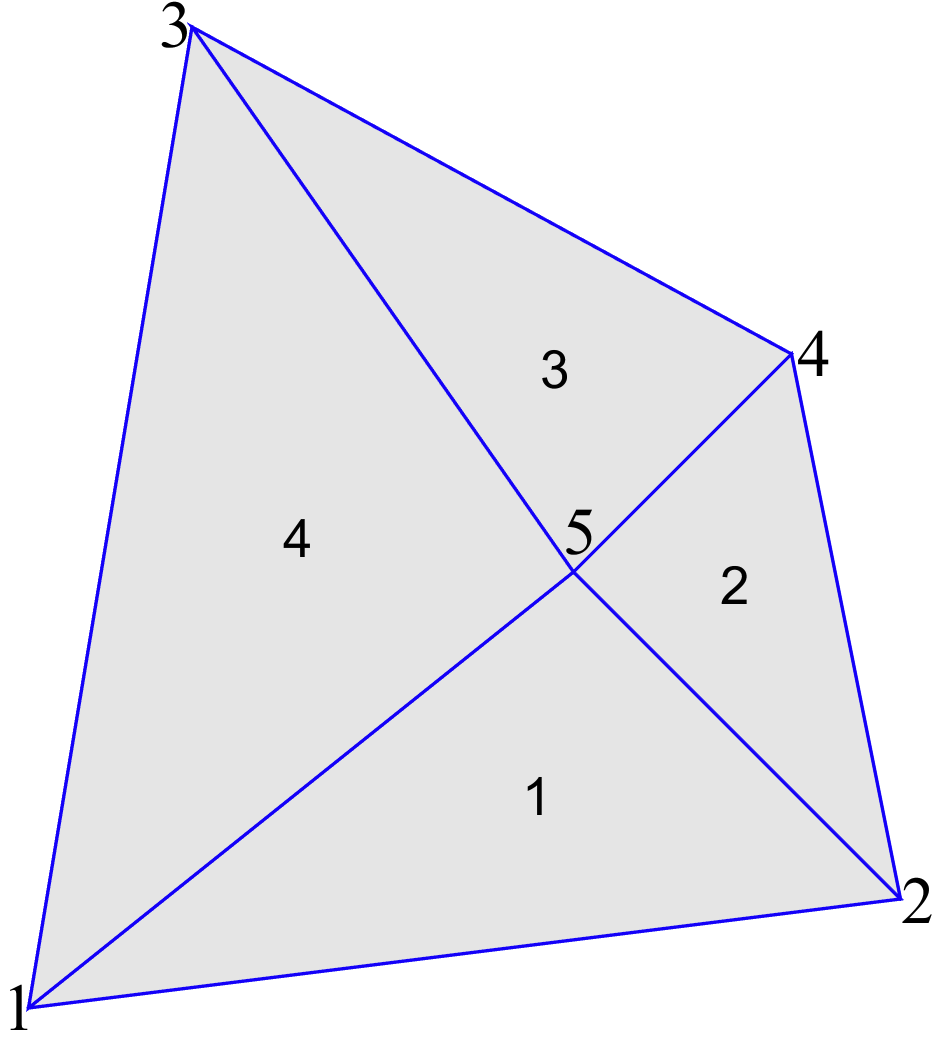}
\caption{Domain discretized with 6 node triangular element}
\label{quadrilateral_mesh_t6_nodes}
\end{subfigure}%
\begin{subfigure}{0.5\columnwidth}
\centering
\includegraphics[width=0.5\textwidth]{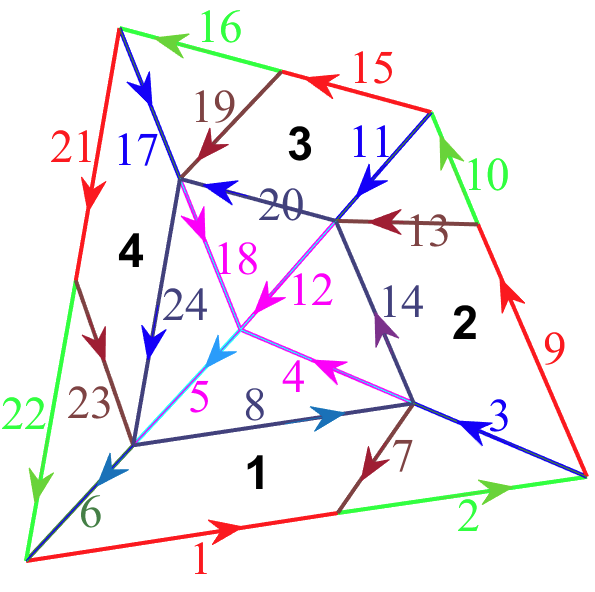}
\caption{Domain discretized with 8 edge triangular element}
\label{quadrilateral_mesh_t6_edges}
\end{subfigure}
\caption{Discretized quadrilateral domain}
\end{figure}

	Higher order triangular edge element (eight edge triangle) is formed from six node triangular element. Fig.~\ref{6_nod_tri} is the local nodal connectivity of six node triangular element. Fig.~\ref{8_edge_tri} represents the desired local edge connectivity of eight edge triangular element. Here eight edges \textbf{e}$_{1}$, \textbf{e}$_{2}$, \textbf{e}$_{3}$, \textbf{e}$_{4}$, \textbf{e}$_{5}$, \textbf{e}$_{6}$, \textbf{e}$_{7}$ and \textbf{e}$_{8}$ are formed by using the node sets \textbf{(1,4)}, \textbf{(4,2)}, \textbf{(2,5)}, \textbf{(5,3)}, \textbf{(3,6)}, \textbf{(6,1)}, \textbf{(6,5)} and \textbf{(5,4)} respectively. $l_1$, $l_2$, ..., $l_8$ are the edge lengths of the element. Here, three edges can be formed on the face of the element. But one edge can be ignored due to independancy of three edges. The edge shape functions of the element are~\cite{High1995} $\bv_1=l_1(4\xi-1)(\xi\del\eta-\eta\del\xi)$, $\bv_2=l_2(4\eta-1)(\xi\del\eta-\eta\del\xi)$, $\bv_3=l_3(4\eta-1)(\eta\del\alpha-\alpha\del\eta)$, $\bv_4=l_4(4\alpha-1)(\eta\del\alpha-\alpha\del\eta)$, $\bv_5=l_5(4\alpha-1)(\alpha\del\xi-\xi\del\alpha)$, $\bv_6=l_6(4\xi-1)(\alpha\del\xi-\xi\del\alpha)$, $\bv_7=4l_7\eta(\alpha\del\xi-\xi\del\alpha)$ and $\bv_8=4l_8\xi(\eta\del\alpha-\alpha\del\eta)$ where $\alpha=1-\xi-\eta$ and $\del\bcross\bE$ can be calculated by using the relation,
\begin{align*}
\del\bcross\bE & =
\begin{bmatrix}
\parder{v_{1y}}{x}-\parder{v_{1x}}{y} & \parder{v_{2y}}{x}-\parder{v_{2x}}{y} & \hdots & \parder{v_{7y}}{x}-\parder{v_{7x}}{y} & \parder{v_{8y}}{x}-\parder{v_{8x}}{y}
\end{bmatrix}
\begin{Bmatrix}
E_1 \\[2mm] E_2 \\[2mm] \vdots \\[2mm] E_{7} \\[2mm] E_{8}
\end{Bmatrix}
= \bB\bE
\end{align*}
where $E_1$, $E_2$, ..., $E_8$ are tangential components of electric fields along the edges \textbf{e}$_1$, \textbf{e}$_2$, ..., \textbf{e}$_8$ respectively. Here, we can obtain the components of the $\bB$-matrix explicitly. We have used the conversion algorithm as discussed in section~\ref{generating_edge_connectivity_array} to convert one general domain as shown in Fig.~\ref{quadrilateral_mesh_t6_nodes} discretized with four 6 node triangular elements into the domain meshed with four eight edge elements shown in Fig.~\ref{quadrilateral_mesh_t6_edges}.    

\subsection{Twelve edge quadrilateral element}

\begin{figure}[pos=h!]
\centering
\begin{subfigure}{0.4\textwidth}
\centering
\includegraphics[width=0.7\textwidth]{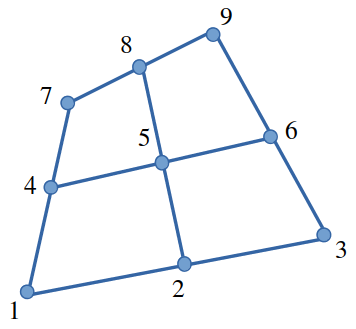}
\caption{Quadrilateral element with nine nodes}
\label{9_nod_quad}
\end{subfigure}%
\begin{subfigure}{0.4\columnwidth}
\centering
\includegraphics[width=0.7\textwidth]{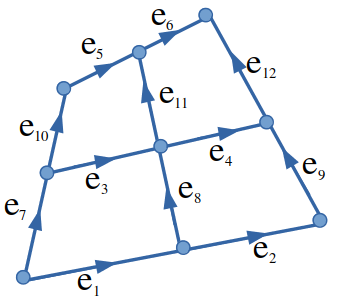}
\caption{Quadrilateral element with twelve edges}
\label{12_edge_quad}
\end{subfigure}
\caption{Higher order quadrilateral element}
\end{figure}
Higher order quadrilateral edge element with twelve edges is formed from the nine node quadrilateral element. The edges \textbf{e}$_{1}$, \textbf{e}$_{2}$, \textbf{e}$_{3}$, ..., \textbf{e}$_{12}$ are formed from the local node sets \textbf{(1,2)}, \textbf{(2,3)}, \textbf{(4,5)}, \textbf{(5,6)}, \textbf{(7,8)}, \textbf{(8,9)}, \textbf{(1,4)}, \textbf{(2,5)}, \textbf{(3,6)}, \textbf{(4,7)}, \textbf{(5,8)} and \textbf{(6,9)} respectively. Local nodal connectivity and edge connectivity are shown in Fig.~\ref{9_nod_quad} and Fig.~\ref{12_edge_quad} respectively. The twelve edge shape functions are~\cite{Jinming-finite} $\bv_1=\frac{-l_1}{2}\eta(\eta-1)(\xi-0.5)\del\xi$, $\bv_2=\frac{l_2}{2}\eta(\eta-1)(\xi+0.5)\del\xi$, $\bv_3=l_3(\eta^2-1)(\xi-0.5)\del\xi$, $\bv_4=-l_4(\eta^2-1)(\xi+0.5)\del\xi$, $\bv_5=\frac{-l_5}{2}\eta(\eta+1)(\xi-0.5)\del\xi$, $\bv_6=\frac{l_6}{2}\eta(\eta+1)(\xi+0.5)\del\xi$, $\bv_7=\frac{-l_7}{2}\xi(\xi-1)(\eta-0.5)\del\eta$, $\bv_8=l_8(\xi^2-1)(\eta-0.5)\del\eta$, $\bv_9=\frac{-l_9}{2}\xi(\xi+1)(\eta-0.5)\del\eta$, $\bv_{10}=\frac{l_{10}}{2}\xi(\xi-1)(\eta+0.5)\del\eta$, $\bv_{11}=-l_{11}(\xi^2-1)(\eta+0.5)\del\eta$ and $\bv_{12}=\frac{l_{12}}{2}\xi(\xi+1)(\eta+0.5)\del\eta$ where $l_1$, $l_2$, ... $l_{12}$ are the lengths of the edges of the element. 

\begin{figure}[pos=h!]
\centering
\begin{subfigure}{0.5\textwidth}
\centering
\includegraphics[width=0.5\textwidth]{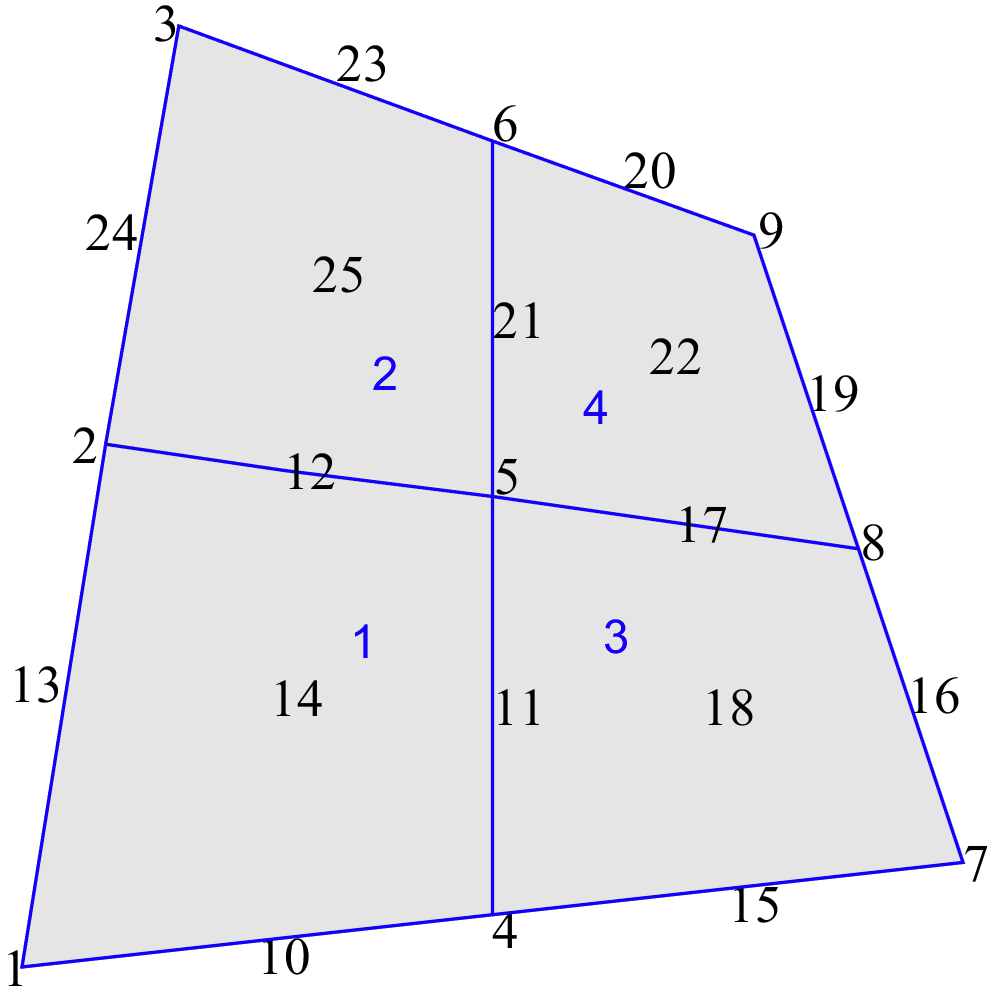}
\caption{Domain discretized with 9 node quadrilateral element}
\label{quadrilateral_mesh_q9_nodes}
\end{subfigure}%
\begin{subfigure}{0.5\columnwidth}
\centering
\includegraphics[width=0.5\textwidth]{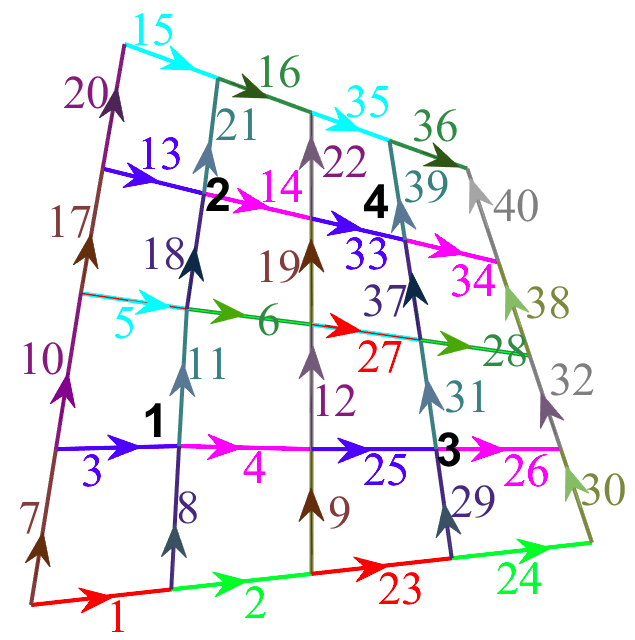}
\caption{Domain discretized with 12 edge quadrilateral element}
\label{quadrilateral_mesh_q9_edges}
\end{subfigure}
\caption{Discretized domain}
\end{figure}
$\del\bcross\bE$ of the element can be calculated by using the following relation.
\begin{align*}
\del\bcross\bE & =
\begin{bmatrix}
\parder{v_{1y}}{x}-\parder{v_{1x}}{y} & \parder{v_{2y}}{x}-\parder{v_{2x}}{y} & \hdots & \parder{v_{11y}}{x}-\parder{v_{11x}}{y} & \parder{v_{12y}}{x}-\parder{v_{12x}}{y}
\end{bmatrix}
\begin{Bmatrix}
E_1 \\[2mm] E_2 \\[2mm] \vdots \\[2mm] E_{11} \\[2mm] E_{12}
\end{Bmatrix}
= \bB\bE
\end{align*}
where $E_1$, $E_2$, ..., $E_{12}$ are tangential components of electric fields along the edges \textbf{e}$_1$, \textbf{e}$_2$, ..., \textbf{e}$_{12}$ respectively. We can obtain the components of the $\bB$-matrix explicitly by using the finite difference method or mathematica software tool~\cite{mathematica}. For the general domain shown in Fig.~\ref{quadrilateral_mesh_q9_nodes} which shows the domain meshed with four 9-node quadrilateral elements we have used the conversion algorithm to transform into the domain as shown in Fig.~\ref{quadrilateral_mesh_q9_edges}. This generated domain is discretized with four twelve edge quadrilateral elements. 

\section{Numerical Examples}
\label{numerical_examples}
In the frequency domain, electromagnetic wave equation can be written as \cite{Jinming-finite}
\begin{equation} \label{eqmaxwell10}
\del\bcross\left (\frac{1}{\mu_r}\del\bcross\bE\right)-k_0^2\epsilon_r\bE=-i\omega\mu_0 \bj,
\end{equation}
where $i=\sqrt{-1}$ and $k_0=\omega/c$ is the wave number in vacuum. From the relations, relative permittivity and relative permeability $\epsilon_r:=\epsilon/\epsilon_0$ and $\mu_r:=\mu/\mu_0$, where $\epsilon_0$ and $\mu_0$ are the permittivity and permeability for vacuum, 
 $c=1/\sqrt{\epsilon_0\mu_0}$ is the speed of light.
Assuming current density, $\bj$ to be zero, the above equation reduces to
\begin{equation} \label{eqmaxwell11}
\del\bcross\left (\frac{1}{\mu_r}\del\bcross\bE\right)=k_0^2\epsilon_r\bE.
\end{equation}
 Eq.~\ref{eqmaxwell11} is utilized to solve the eigenvalue problems for finding the square of the eigenvalue $k_0^2$. In order to validate the transformed edge elements we have performed numerical analysis by considering the standard eigenvalue problems. Here, edge elements such as 3-edge triangular, 4-edge quadrilateral, 12-edge quadrilateral and 8-edge triangular elements are represented as T3, Q4, Q12 and T8 respectively. We assume $\epsilon_r=\mu_r=1.0$ for all the problems considered under this section.

\subsection{Square domain with perfectly conducting boundaries}    

Square domain of side length $\pi$ is chosen. For this square domain all the sides/boundaries are assumed to be perfectly conducting. Here, to conduct the numerical analysis, the domain is discretized with the first order edge elements (T3 and Q4) and higher order edge elements (Q12 and T8). Table~\ref{square_analysis_data} shows the total number of equations along with the total number of discretized elements for different meshes. For each element square of eigen values are listed in Table~\ref{tabsquare1}. Results with the transformed edge elements are in good agreement with the analytical results stated in ~\cite{benchmarkmonique}. We can observe that all the elements gave right multiplicity of eigenvalues. The first non-zero eigenvalues of all the elements appeared after stating the number of zeros generated at the machine precision level. These zeros indicate the approximation of null space.  
\begin{table}[width=.79\textwidth,pos=ht!]
\caption{Analysis data of different edge elements for the square domain problem.} \label{square_analysis_data} 
\begin{tabular*}{\tblwidth}{@{}LLL@{}} \toprule  
\textbf{Type of} & \textbf{Total no. of} & \textbf{Total free degrees}  \\ 
\textbf{Element} & \textbf{Elements} & \textbf{of freedom (equations)} \\ \midrule  
T3  & 512  & 736  \\  
Q4  & 256  & 480  \\  
T8  & 128 & 1280  \\  
Q12 & 256  & 1984  \\ \bottomrule 
\end{tabular*}
\end{table}

\begin{table}
\caption{$k_0^2$ on the square domain for different elements.} \label{tabsquare1}
\begin{tabular*}{\tblwidth}{@{}LLLLL@{}}  \toprule
\textbf{Analytical} & \multicolumn{4}{c}{\textbf{Edge element}}  \\ \midrule 
\textbf{Benchmark} & \textbf{T3} & \textbf{Q4} & \textbf{Q12} & \textbf{T8} \\ \midrule  
 1   &     0.998066  &    1.000803  &    1.000002   &   0.999992 \\  
 1   &     0.999795  &    1.000803  &    1.000002   &   1.000010 \\  
 2   &     2.002121  &    2.001607  &    2.000004   &   2.000115 \\  
 4   &     3.982881  &    4.012868  &    4.000131   &   4.000089 \\  
 4   &     3.982939  &    4.012868  &    4.000131   &   4.000089 \\  
 5   &     4.982602  &    5.013671  &    5.000133   &   5.000260 \\  
 5   &     5.015107  &    5.013671  &    5.000133   &   5.002108 \\  
 8   &     8.032183  &    8.025735  &    8.000262   &   8.006889 \\  
 9   &     8.906076  &    9.065245  &    9.001478   &   9.000147 \\  
 9   &     8.921107  &    9.065245  &    9.001478   &   9.001707 \\  
 10   &      9.950139  &   10.066048  &   10.001480   &  10.005688 \\  
 10   &      9.952486  &   10.066048  &   10.001480   &  10.005711 \\  
 13   &     12.960172  &   13.078112  &   13.001609  &   13.012005  \\
 13   &     13.133842  &   13.078112  &   13.001609  &   13.037121  \\
 16   &     15.726881  &   16.206657  &   16.008194  &   16.004350  \\
 16   &     15.727173  &   16.206657  &   16.008194  &   16.004383  \\  \midrule
\multicolumn{2}{c}{Number of computed zeros}  &  \multicolumn{3}{c}{}       \\  \midrule
                   -        &       65       &     220      &   217        &    224       \\ \bottomrule
\end{tabular*}
\end{table}

\subsection{Curved-L shape domain with perfectly conducting surfaces}    

\begin{figure}[h!]
\centering
\includegraphics[width = 12.0 cm, height = 6.0 cm]{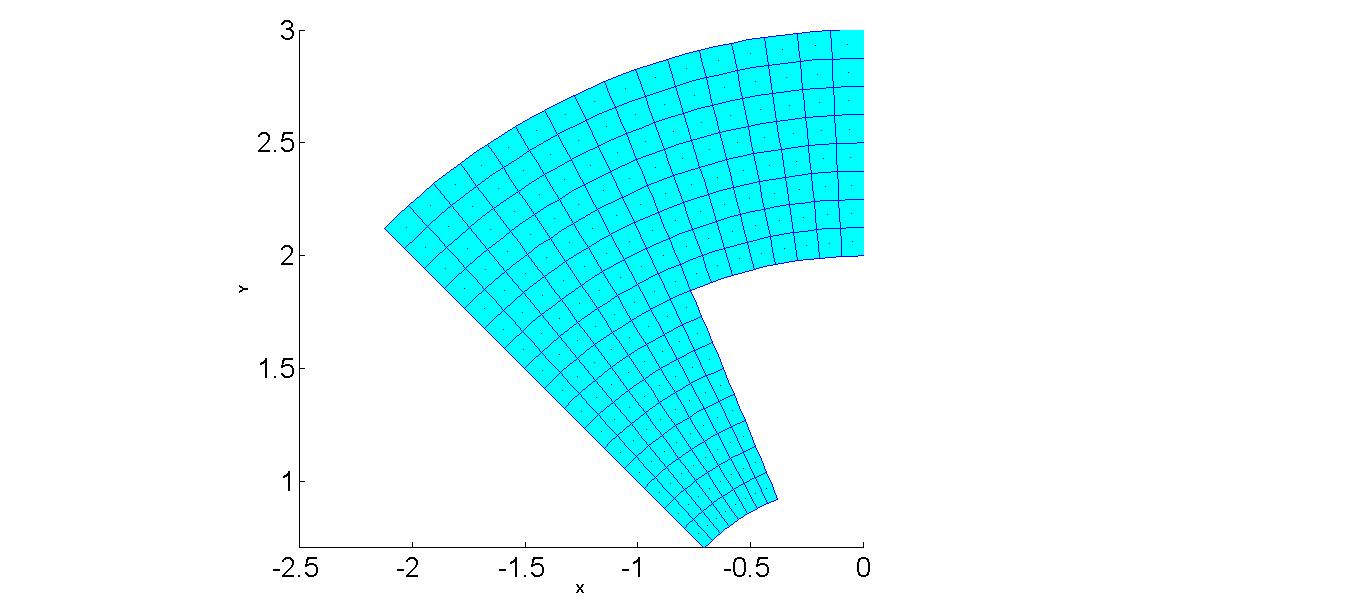}
\caption{Mesh for Curved L-shape domain}
\label{lshape}
\end{figure}
 In this example, solved in~\cite{benchmarkmonique}, the domain has three straight and three circular sides of radii 1, 2 and 3 and Fig.~\ref{lshape} shows such domain discretized with Q4 elements. Here, all the boundary edges of the domain are perfectly conducting. This problem is quite complicated due to the existence of singular eigen value for the sharp corner, and curvature effect. To find the eigen values we discretize the domain with the T3, Q4, Q12 and T8 elements. In Table~\ref{curvedL_data} we presented total number of equations and total number of elements for each type of element. For all the elements numerical results obtained are listed in Table~\ref{tablshape2} along with the number of computed zeros. We compared the numerical results of these elements with analytical values taken from ~\cite{benchmarkmonique}. It can be observed that higher order edge element Q12 results matches with the benchmark values up to second decimal and for T8 elements it matches upto first decimal.   

\begin{table}[width=.79\textwidth,pos=ht!]
\caption{Analysis data of different edge elements for the curved-L shape domain problem.} \label{curvedL_data} 
\begin{tabular*}{\tblwidth}{@{}LLL@{}} \toprule  
\textbf{Type of} & \textbf{Total no. of} & \textbf{Total free degrees}  \\ 
\textbf{Element} & \textbf{Elements} & \textbf{of freedom (equations)} \\ \midrule  
T3  & 600   & 860  \\  
Q4  & 300   & 560  \\  
T8  & 96    & 448  \\  
Q12  & 108  & 816  \\ \bottomrule 
\end{tabular*}
\end{table}

\begin{table}[width=.79\textwidth,pos=ht!]
\caption{$k_0^2$ on the curved L-shaped domain for different elements.} \label{tablshape2}
\begin{tabular*}{\tblwidth}{@{}LLLLL@{}}  \toprule
\textbf{Analytical} & & \textbf{Edge element} & & \\ \midrule 
\textbf{Benchmark} & \textbf{T3} & \textbf{Q4} & \textbf{Q12} & \textbf{T8} \\ \midrule  
  1.818571   &  1.797075   &  1.811631   & 1.814860   &  1.807729   \\
  3.490576   &  3.491215   &  3.500850   & 3.490516   &  3.4954251  \\
 10.065602   &  10.047041  &  10.151037  & 10.066760  &  10.082410  \\
 10.111886   &  10.101835  &  10.203259  & 10.112480  &  10.127115  \\
 12.435537   &  12.397735  &  12.510154  & 12.429986  &  12.431268  \\ \midrule
\multicolumn{2}{c}{Number of computed zeros}  &  \multicolumn{3}{c}{}       \\  \midrule
                              -      &   101       &     220     &    217     &    224      \\ \bottomrule
\end{tabular*}
\end{table}


\subsection{Circular domain with perfectly conducting surfaces}
\label{circular_domain}

\begin{figure}[h!]
\centering
\includegraphics[width = 12.0 cm, height= 6.0 cm]{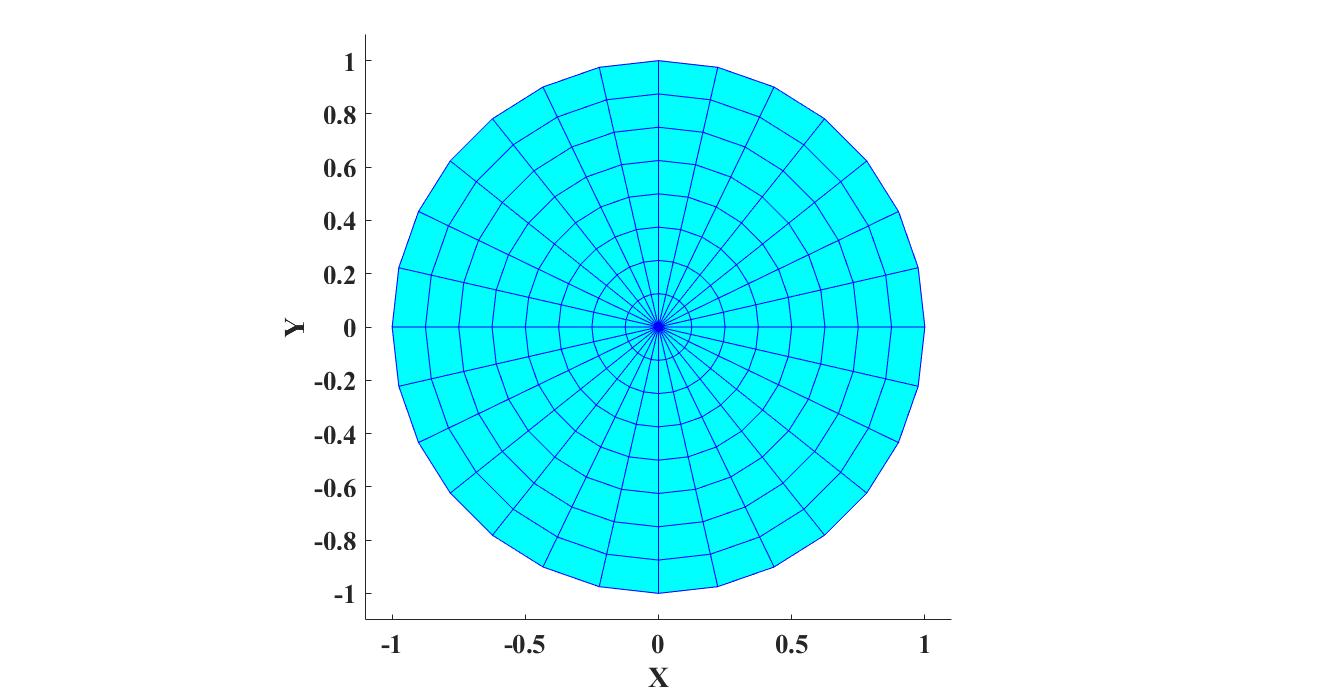}
\caption{Discretized circular domain}
\label{cir_dom}
\end{figure}

A circular domain of unit radius with perfectly conducting boundaries is considered to perform eigen analysis. Our interest to consider this domain is to test the proposed algorithm in handling the data of the additional edges (more than four edges shared at a particular node) as mentioned in section~\ref{nodeedgeexn_circle_domain}. In the earlier examples the whole domains are discretized with only one type of elements. But in the present case the domain is discretized with the combination of first order edge elements (T3 and Q4) or the combination of higher order edge elements (Q12 and T8) to perform numerical analysis. Triangular elements are used to mesh the centre portion of the domain upto one layer in $r$ direction. Quadrilateral elements are adopted to discretize the rest of the domain as shown in Fig.~\ref{cir_dom}. Mesh details of these elements are shown in Table~\ref{circular_data} and $k_0^2$ values are listed in Table~\ref{tabcir} and are compared with analytical results reported in ~\cite{2554},~\cite{harrington2001time},~\cite{Jog2014}. The results of the generated edge elements indicate a strong fit with the analytical results along with the correct multiplicity of eigen values.

\begin{table}[width=.79\textwidth,pos=ht!]
\caption{Analysis data of different edge elements for the circular domain problem.} \label{circular_data} 
\begin{tabular*}{\tblwidth}{@{}LLLL@{}} \toprule  
\textbf{Type of} & \textbf{Total no. of} & \textbf{Total no. of} & \textbf{Total free degrees}  \\ 
\textbf{Element} & \textbf{Triangular elements} & \textbf{Quadrilateral elements} & \textbf{of freedom (equations)} \\ \midrule  
Q4/T3  & 30 & 1770  & 1430  \\  
Q12/T8 & 20 & 680  & 6340  \\ \midrule 
\end{tabular*}
\end{table}

\begin{table}
\centering
\caption{$k_0^2$ on the circular domain for different elements} \label{tabcir}
\begin{tabular*}{\tblwidth}{@{}LLL@{}}  \toprule
\textbf{Analytical} & \multicolumn{2}{c}{\textbf{Edge element}}  \\ \midrule 
\textbf{Benchmark} & \textbf{T3/Q4} & \textbf{T8/Q12}  \\ \midrule  
3.391122(2)    & 3.425827(2)  & 3.383070(2)  \\ 
9.329970(2)    & 9.530267(2)  & 9.329383(2)  \\ 
14.680392(1)   & 14.802788(1) & 14.737119(1) \\ 
17.652602(2)   & 18.347489(2) & 17.662205(2) \\ 
28.275806(2)   & 28.658911(2) & 28.304885(2) \\ 
28.419561(2)   & 30.105560(2) & 28.343554(2) \\ 
41.158640(2)   & 45.175889(2) & 41.404575(2) \\ 
44.970436(2)   & 45.564095(2) & 44.974232(2) \\ 
49.224256(1)   & 49.668776(1) & 49.543520(1) \\ 
56.272502(2)   & 64.056343(2) & 56.956832(2) \\ 
64.240225(2)   & 65.762272(2) & 64.271185(2) \\  \midrule 
\multicolumn{2}{c}{Number of computed zeros}  &  \multicolumn{1}{c}{}       \\  \bottomrule
                              -      &       494     &     629  \\ \hline
\end{tabular*}
\end{table}

\subsection{Cracked circular domain with perfectly conducting surfaces}
\label{cracked_circular_domain}

\begin{figure}[h!]
\centering
\includegraphics[width = 5.0 cm, height= 5.0 cm]{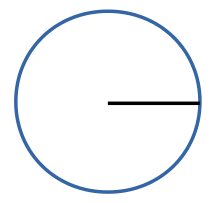}
\caption{Cracked circular domain}
\label{crack_cir_dom}
\end{figure}

In this example, same circular domain but with the crack running from the centre to the side of the circle is taken into consideration, as shown in Fig.~\ref{crack_cir_dom}. Here, the domain's crack is modeled by using `double noding' method at same position. Mesh details are shown in Table~\ref{crack_circle_data} including the total number of equations. Table~\ref{tabcrackcir} shows the numerical values along with the computed zeros and they are compared with analytical results reported in ~\cite{2554},~\cite{harrington2001time},~\cite{Jog2014}. The results of the edge elements (T8/Q12) demonstrate close matching with the analytical results.

\begin{table}[width=.79\textwidth,pos=h!]
\caption{Analysis data of different edge elements for the cracked circular domain problem.} \label{crack_circle_data} 
\begin{tabular*}{\tblwidth}{@{}LLLL@{}} \toprule  
\textbf{Type of} & \textbf{Total no. of} & \textbf{Total no. of} & \textbf{Total free degrees}  \\ 
\textbf{Element} & \textbf{Triangular elements} & \textbf{Quadrilateral elements} & \textbf{of freedom (equations)} \\ \midrule  
Q4/T3  & 30 & 1770  & 3631  \\  
Q12/T8 & 20 & 680  & 5591  \\ \midrule 
\end{tabular*}
\end{table}

\begin{table}
\centering
\caption{$k_0^2$ on the cracked circular domain for different elements} \label{tabcrackcir}
\begin{tabular*}{\tblwidth}{@{}LLL@{}}  \toprule
\textbf{Analytical} & \multicolumn{2}{c}{\textbf{Edge element}}  \\ \midrule 
\textbf{Benchmark} & \textbf{T3/Q4} & \textbf{T8/Q12}  \\ \midrule  
1.358390    & 1.362745  & 1.297322 \\ 
3.391122    & 3.425901  & 3.383104 \\ 
6.059858    & 6.146835  & 6.053942 \\ 
9.329970    & 9.530267  & 9.329383 \\ 
13.195056    & 13.589143 &  13.201275\\ 
14.680392    &    -      &  15.304855\\ 
17.652602    & 18.347489 &  17.662205\\ 
21.196816    & 21.297238 &  20.547667\\ 
22.681406    & 23.838848 &  22.709137\\ 
28.275806    & 28.660201 &  28.305522\\ \midrule
\multicolumn{2}{c}{Number of computed zeros}  &  \multicolumn{1}{c}{}       \\  \midrule
                              -      & 546     &     659  \\ \bottomrule
\end{tabular*}
\end{table}

\section{Conclusions}
Electromagnetic analysis with nodal finite element has several shortcomings. Nodal FEM cannot model the null space accurately, there are presence of spurious values. With regularization or penalty method this spurious values are shifted towards higher end. But with this method an adhoc penalty parameter is required to be adjusted. Also, with this penalty method singular eigen values for the domains with sharp edges and corners, cannot be approximated accurately. With edge finite element method all these limitations are addressed without the use of any adhoc penalty parameter. Furthermore, in nodal FEM, it is required to decompose electric and magnetic fields into scalar and vector potentials to attain the necessary continuity requirement across elements. After FEM analysis, fields are calculated from the potentials with additional postprocessing. In edge FEM, we can formulate directly in terms of field variables. But most of the preprocessor in practice, generate FEM meshes in terms of nodal connectivities. Hence, in this article we have presented a very useful novel conversion technique which transform the nodal connectivities into edge connectivities. Also, this conversion algorithm generates other necessary data structures in edge formulations like direction information of the edges, connecting nodes of a particular edges, associated edges of a particular node and respective other node of those edges. This algorithm converts 4-node quadrilateral into 4-edge quadrilateral, 3-node triangle into 3-edge triangle, 6-node triangle into 8-edge triangle, and 9-node quadrilateral into 12-edge quadrilateral. For some special geometries combination of triangular and quadrilateral elements are more effective. Our conversion algorithm is capable to combine successfully 3-edge triangle with 4-edge quadrilateral and 8-edge traingle with 12-edge quadrilateral. In section~\ref{circular_domain}  and~\ref{cracked_circular_domain} successful implementation of such combination is presented with numerical examples. Some special treatment is required in the data structure for the nodes which are connected to many edges; it is explained in detail in section~\ref{nodeedgeexn_circle_domain} with associated representative examples. This additional data structure is verified with numerical examples in section~\ref{circular_domain} and~\ref{cracked_circular_domain}. The effectiveness of the conversion technique is tested with different standard benchmark examples. These numerical examples include square domain, circular domain, cracked circular domain, curved L shape domain etc. Our proposed conversion technique gives accurate $k_0^2$ values along with correct multiplicity for both convex and non-convex domains. For non-convex domains, singular eigenvalues are predicted without any spurious modes. The perfect match with the benchmark results for different examples in terms of eigenvalues, their multiplicities, singular eigenvalues for domain with sharp corners and edges, exhibit the correctness and efficacies of the conversion algorithm. 

\bibliographystyle{cas-model2-names}

\bibliography{library.bib}

\end{document}